\documentclass[sn-mathphys]{sn-jnl}

\usepackage{longtable}
\usepackage{verbatim}
\usepackage{cancel}

\jyear{2022}%

\theoremstyle{thmstyleone}%
\newtheorem{theorem}{Theorem}
\newtheorem{proposition}[theorem]{Proposition}%

\theoremstyle{thmstyletwo}%
\newtheorem{remark}{Remark}%

\theoremstyle{thmstylethree}%
\newtheorem{definition}{Definition}%
\newtheorem{assumption}{Assumption}
\newtheorem{lemma}{Lemma}
\newenvironment{Proof}{\noindent {\it Proof. }}{\hfill $\square$}

\newcommand{\cX}{{\cal X}}
\newcommand{\cT}{{\cal T}}
\newcommand{\cJ}{{\cal J}}
\newcommand{\cF}{{\cal F}}
\newcommand{\cG}{{\cal G}}
\newcommand{\cU}{{\cal U}}

\newcommand{\cM}{{\cal M}}
\newcommand{\cN}{{\cal N}}

\newcommand{\cA}{{\cal A}}

\newcommand{\cQ}{{\cal Q}}
\newcommand{\cB}{{\cal B}}

\newcommand{\bea}{\begin{eqnarray*}}
	\newcommand{\eea}{\end{eqnarray*}}
\newcommand{\mb}[1]{\mbox{\boldmath $#1$}}

\newcommand{\inprod}[2]{\langle #1 , #2 \rangle }

\newcommand{\bc}{\begin{center}}
	\newcommand{\ec}{\end{center}}
\newcommand{\bz}{\vect{z}}

\newcommand{\bv}{\vect{v}}
\newcommand{\by}{\vect{y}}

\newcommand{\bx}{\vect{x}}

\newcommand{\RR}{\mathbb R}

\newcommand{\be}{\begin{equation}}
\newcommand{\ee}{\end{equation}}
\newcommand{\beaa}{\begin{eqnarray*}}
	\newcommand{\eeaa}{\end{eqnarray*}}
\newcommand{\ben}{\begin{enumerate}}
	\newcommand{\een}{\end{enumerate}}
\newcommand{\db}{\hspace*{\fill}{\zapf o}}

\newcommand{\bpn}{\begin{proposition}\twlsf}
	\newcommand{\epn}{\db\end{proposition}}
\newcommand{\bdm}{\begin{displaymath}}
\newcommand{\edm}{\end{displaymath}}
\newcommand{\ba}{\begin{array}}
	\newcommand{\ea}{\end{array}}

\newcommand{\argmin}{\mathop{\rm argmin}}

\def\texitem#1{\par\smallskip\noindent\hangindent 25pt
	\hbox to 25pt {\hss #1 ~}\ignorespaces}
\newcommand{\norm}[1]{\left\lVert#1\right\rVert}

\newcommand{\eps}{\epsilon}

\newcommand{\blue}[1]{\begin{color}{black}#1\end{color}}

\def\lam{\lambda} \def\alp{\alpha}
\def\sig{\sigma}

\def\nn{\nonumber}

\def\inprod#1#2{\langle#1,\,#2\rangle}
\def\norm#1{\|#1\|}
\def\diag#1{\mbox{diag}(#1)}
\def\grad{\nabla}

\def\S{\mathbb{S}} 
\def\cL{{\cal L}}  \def\cK{{\cal K}}  \def\cU{{\cal U}}
\def\cY{{\cal Y}}   \def\cD{{\cal D}}
\def\cB{{\cal B}}

\def\bv{\bar{v}}
\def\bx{\bar{x}} \def\by{\bar{y}} \def\bz{\bar{z}}
 \def\tx{\tilde{x}} 
\def\tc{\tilde{c}} \def\tB{\tilde{B}}
\def\tb{\tilde{b}}
\def\ctX{\tilde{\cX}} \def\ctY{\tilde{\cY}}

\def\bx{\bar{x}} \def\bz{\bar{z}} 
   
\def\btheta{\bar{\theta}}
\def\bcK{{\bar{\cal K}}}
\def\cX{{\cal X}} \def\bcX{{\bar{\cal X}}}  \def\bcY{{\bar{\cal Y}}} 
\def\cU{{\cal U}} 
 \def\cN{{\cal N}} \def\cL{{\cal L}}
\def\cS{{\cal S}}
\def\bb{\bar{b}} \def\bc{\bar{c}}\def\bB{\bar{B}} \def\bQ{{\bar{Q}}}
 
\def\bm{\bar{m}} \def\bn{\bar{n}} \def\by{\bar{y}}  
 \def\bD{\bar{D}}  
 \def\bJ{\bar{J}}
\def\bE{\bar{E}}

\def\mbU{\mb{U}}

\begin{document}

\title[Proximal ALM for dual block-angular convex composite  programming]{On proximal augmented Lagrangian based decomposition methods for dual block-angular convex composite  programming problems}


\author[1]{\fnm{Kuang-Yu} \sur{Ding}}\email{kuangyud@u.nus.edu}

\author*[1]{\fnm{Xin-Yee} \sur{Lam}}\email{matlxy@nus.edu.sg}
\equalcont{These authors contributed equally to this work.}

\author[1,2]{\fnm{Kim-Chuan} \sur{Toh}}\email{mattohkc@nus.edu.sg}
\equalcont{These authors contributed equally to this work.}

\affil*[1]{\orgdiv{Department of Mathematics}, \orgname{National
		University of Singapore}, \orgaddress{\street{10 Lower Kent Ridge Road}, \postcode{119076}, \country{Singapore}}}

\affil[2]{\orgdiv{Institute of 
		Operations Research and Analytics}, \orgname{National
		University of Singapore}, \orgaddress{\street{3 Research Link}, \postcode{117602}, \country{Singapore}}}


\abstract{We design {\em inexact} proximal augmented Lagrangian based decomposition methods for 
	convex composite  programming problems with dual block-angular structures. 
	Our methods are particularly well suited for convex quadratic  programming problems
	arising from stochastic programming models. 
	The algorithmic framework is based on {the application of
		the abstract inexact proximal ADMM framework developed in [Chen, Sun, Toh, Math. Prog. 161:237--270] to the dual of the target problem, as well as the application of the recently
		developed
		symmetric Gauss-Seidel decomposition theorem for solving a proximal multi-block convex composite quadratic programming problem.
		The key issues in our algorithmic design are firstly in designing appropriate 
		proximal terms to  
		decompose the computation of the dual variable blocks of the target problem
		to make the subproblems in each iteration easier to solve, and secondly to develop 
		novel numerical schemes to solve the decomposed subproblems efficiently.}
	Our inexact augmented Lagrangian based decomposition methods have guaranteed convergence. We present an application of the proposed algorithms to the doubly nonnegative relaxations of uncapacitated facility location problems, as well as to two-stage stochastic optimization problems. We conduct numerous numerical experiments to 
	evaluate the performance of our method against state-of-the-art solvers such as Gurobi and MOSEK. Moreover, our proposed algorithms also compare favourably 
	to the well-known progressive hedging algorithm of Rockafellar and Wets. }




\maketitle

\newpage 
\section{Introduction}\label{sec1}

Let $\mathcal X$ and $\cX_i$ ($i=1,\ldots,N$) 
be $n_0$-dimensional and $n_i$-dimensional  inner product spaces each endowed with the inner product
$\inprod{\cdot}{\cdot}$ and its induced norm $\norm{\cdot}$, respectively. 
Similarly, let $\cY$ and $\cY_i$ ($i=1,\ldots,N$) be $m_0$-dimensional and $m_i$-dimensional inner product spaces,
respectively.
We consider the following general convex composite  programming  problem with the so called
dual block-angular structure (DBA) (according to the terminology used by Dantzig \cite{Dantzig} in the case of linear programming): 
\begin{eqnarray}
\begin{array}{rl}
\min & \theta(x) + \inprod{c}{x} + \sum_{i=1}^N \blue{\left(\btheta_i(\bx_i)+\inprod{\bc_i}{\bx_i}\right) } \\[5pt]
\mbox{s.t.}  & \left[ \begin{array}{ccccc}
A &  \\[5pt]
B_1 &  \bB_1  &  \\
B_2 &   & \bB_2 &   \\
\vdots  &     &  &\ddots\\
B_N     &     &      & & \bB_N
\end{array} \right]
\left[ \begin{array}{c}
x \\ \bx_1 \\ \vdots \\ \bx_N
\end{array} \right]
\;=\;
\left[ \begin{array}{c}
b \\ \bb_1 \\ \vdots \\ \bb_N
\end{array} \right]
\\[40pt] 
& x \in \cK, \quad  \bx_i \in \cK_i \subset \cX_i, \;\forall\; i=1,\ldots,N,
\end{array}
\label{eq-DBA}
\end{eqnarray}
where $b\in \cY$, $c\in\cX$, $\bb_i\in \cY_i$ and $c_i \in \cX_i$ ($i=1,\ldots,N$) 
are given data; $A:\cX\to \cY$, $B_i:\cX \to\cY_i$ and $\bB_i:\cX_i \to \cY_i$ are 
given linear maps; 
$\theta:\cX\to (-\infty,\infty]$ and $\btheta_i:\cX_i\to (-\infty,\infty]$  are 
given proper closed convex functions; 
$\cK\subset\cX$ and $\cK_i\subset\cX_i$ ($i=1,\ldots, N$)  are given closed convex sets, 
\blue{which are typically cones but not necessarily so.}
In the above problem, we allow the case where $A$ is absent. 
We will call the form in (1) as the standard DBA problem. In the literature, there are
other block-angular forms such as the standard primal block-angular form and other variants
which are not expressed in the standard form. Following the terminology in the literature,
we will sometimes call the variable $x$ as the first-stage variable, and $\bar{x}_i$ ($i = 1,\ldots, N$) as the
second-stage variables.

The above problem (which is reformulated as (P) {in next page})
that we consider in this paper is very general. For example, it includes 
block-angular
convex quadratic programming problems as a special case when 
$\theta(\cdot)$ and $\bar{\theta}_i(\cdot)$  
($i = 1, \ldots,N$) are convex quadratic functions, i.e, 
$\theta(x) = \frac{1}{2}
\inprod{x}{Qx}$ and $\bar{\theta}_i(\bar{x}_i) = \frac{1}{2}\inprod{\bar{x}_i}{\bar{Q}_i\bar{x}_i}$
where $Q : \cX \to \cX$ and $\bar{Q}_i :\cX_i \to \cX_i$
are self-adjoint positive semidefinite linear operators. As
we shall see later, our standard DBA problem (1) is much more general than other forms that
have been considered so far in the literature. 

The dual block-angular problem generally has a separable convex objective function and {convex} constraints but the linear constraints are coupled by the linking variable ($x$) across different constraints. The dimension of these problems can easily grow very large especially when the number of blocks, $N$, increases. In particular, the number of linear equality constraints
is equal to $\sum_{i=0}^N m_i$. 
The main goal of this paper is
to develop a specialized algorithm that can fully exploit the block-angular structure
to solve the large-scale problem efficiently rather than treating it as a generic {convex}
programming problem for which even a powerful commercial solver may not be
able to handle efficiently. 

In designing decomposition based methods, it is important to clearly specify the block-angular
structure one is dealing with, since the decomposition strategies and the efficiency of the
methods are invariably dictated by the structure. In this paper, we will only be concerned with
the DBA structure presented in (1), but not touch on other block-angular structures as they will
lead to different decomposition strategies and methods. Consequently, our literature review will
also be confined to the DBA problems. Even then, the literature on DBA problems is still too
vast for us to give a comprehensive coverage. Thus we shall only touch on papers which are
most directly relevant to our current paper when we described the related literature in Section 1.1.

In this paper, we will use the notation $[P; Q]$ to denote the 
matrix obtained by appending the matrix $Q$ to the last row of the matrix $P$,
assuming that  they have the same number of columns. We will also use the same notation
symbolically
for linear maps $P$ and $Q$ as long as they have the compatible domain and co-domain. 
Let $\bm =\sum_{i=1}^N m_i $  and $\bn= \sum_{i=1}^N n_i $ with $n_i$ being the dimension of $\cX_i$. 
Define 
\begin{eqnarray}
& \bcX =  \cX_1\times \cdots\times \cX_N, \quad \bcY = \cY_1\times \cdots\times \cY_N, \quad 
\bcK = \cK_1\times \cdots\times \cK_N,
&
\nn \\[5pt]
& \bx = [\bx_1;\dots;\bx_N] \in \bcX,
\;\; \bc =[\bc_1;\dots;\bc_N]\in \bcX,\; \;
\bb = [\bb_1;\dots ; \bb_N]\in \bcY, & \nn
\\[5pt]
& B = [B_1;\dots; B_N], \; \;
\bB = \diag{\bB_1,\ldots,\bB_N}, \quad
\btheta(\bx) = \sum_{i=1}^N \btheta_i (\bx_i).
&
\end{eqnarray} 
Note that 
the dimension of $\bcX$ is $\bn$ and that of $\bcY$ is $\bm.$
For convenience, we will sometimes refer to $\bB_i$ as an $m_i\times n_i$ matrix
to mean its matrix representation with respect to the standard bases in $\cX_i$ and $\cY_i$. 
A similar convention may also be used for $A$ and $B_i$. 

We can rewrite \eqref{eq-DBA} as an instance of
the following general  {convex} programming problem with a  block-angular structure:
\begin{eqnarray*}
	\begin{array}{rcll} 
		\mbox{(P)} \quad \min & \theta(x) + \delta_\cK(x) + \inprod{c}{x} 
		&+\; \btheta(\bx)+\delta_{\bcK}(\bx) + \inprod{\bc}{\bx}  \\[3pt]
		\mbox{s.t.} & A x &&= b, \quad x\in \cX  \\[3pt]
		& Bx  &+\; \bB \bx &=\bb, \quad \bx \in \bcX 
	\end{array}
	\label{eq-pp}
\end{eqnarray*}
where $A:\cX \to \cY$, $B:\cX\to \bcY$, $\bB:\bcX \to \bcY$, and $\delta_\cK(\cdot)$, $\delta_\bcK(\cdot)$
denote the indicator functions over $\cK$, $\bcK$, respectively. 
The above structure is often referred to in the literature as a
dual block-angular structure. But we shall drop the word ``dual'' in this paper. 

We shall show in Section \ref{sec:dual} that the dual of (P) is given by
\begin{eqnarray*}
	\begin{array}{rrl}
		\mbox{(D)} \quad   -\min & \theta^*(-v) + \delta_\cK^*(-z) 
		-\inprod{b}{y} \;+ \btheta^*(-\bar{v}) + \delta_{\bcK}^*(-\bz) -\inprod{\bb}{\by}   \\[5pt]
		\mbox{s.t.} & A^* y + B^*\by +z+ v = c, \;\; z\in \cX,\; v\in\cX  \\[5pt]
		&  \;\;\; \bB^*\by+ \bz + \bar{v}  = \bc, \;\; \bz\in \bcX,\; \bar{v}\in \bcX
	\end{array}
	\label{eq-dd}
\end{eqnarray*}
where $ \bz = [\bz_1;\dots;\bz_N] \in \bcX,$ 
$\bar{v} = [\bar{v}_1; \ldots; \bar{v}_N]\in \bcX$, $\bar{y} = [\bar{y}_1;\dots; \bar{y}_N] \in 
\bcY$,  
$\delta^*_\cK$ is the conjugate function of $\delta_\cK$ (similarly for  
$\delta^*_\bcK$). 
Observe that the problem (D) has six blocks of variables $(y,\by,z,\bz,v,\bar{v})$
coupled together by the linear constraints 
$A^*y + B^*\by + z + v = c$ and $\bB^*\by + \bz + \bar{v} = \bc$. 

\bigskip
We made the following assumptions on the problem data throughout this paper. 
\begin{assumption} 
	$A\in \RR^{m\times n}$ has full row rank (if it is present), and $[B,\bar{B}]$ has full row rank.
\end{assumption}
\begin{assumption} 
	The projection $\Pi_\cK(x)$ and $\Pi_\bcK(\bz)$ can be  computed efficiently such as in the case when 
	$\cK$ is a nonnegative orthant, a second-order cone or $\S^n_+$ (the cone of symmetric positive semidefinite 
	matrices). We also assume that the proximal mappings for $\theta$ and $\btheta$ can be computed
	analytically or very efficiently. 
\end{assumption}

In practical applications, optimization problems with dual block-angular structures arise in many contexts {\cite{ZFEM,FLS,Wollenberg}}. In this paper, we focus 
on two classes of problems, namely, two-stage stochastic {convex} programming problems and doubly nonnegative relaxations of uncapacitated facility location problems, for 
the purpose of evaluating the numerical performance of our algorithm against
other existing state-of-the-art solvers. 

In the case of two-stage stochastic {convex} programming problems, the optimization model often has a linear or quadratic objective function. Commonly used methods to solve two-stage stochastic linear-quadratic programming problems can mainly {be} classified 
into three types. The first type is based on the cutting plane method, including Bender decomposition methods \cite{Bender} {and} L-shaped methods \cite{Lshape}. However, this class of methods 
is known to have very slow convergence {\cite{TLAG}}.
The second type is based on  interior point methods, where details can be found in 
\cite{BeDeOZ99, birge95, zhaogy991, zhao2001} and Mehrotra's works \cite{SM1,SM2}. 
The last type is decomposition methods based on augmented Lagrangian method (ALM), see for example \cite{RTRW, CDZ,ruszcz95,ruszcz99} and the references therein.  Others include distributed simplex methods \cite{ParSplx}, and Karmarkar algorithms \cite{Karm,BirQi88}. 
Among the decomposition methods, the Progressive Hedging Algorithm (PHA) in \cite{RTRW}, which is a scenario decomposition based method, is perhaps the best known for 
solving 2-stage stochastic linear-quadratic programming problems.  
For the class of uncapacitated facility location problems, they are originally binary integer programming problems. We relax them to doubly nonnegative (DNN) {convex} programming problems and solve the resulting dual block-angular DNN problems. 

The problem (P) is generally large scale and developing efficient algorithms that are capable of handling such a large problem is of great interest.
In this paper, our main objective is to design and implement an efficient and robust (distributed) algorithm for solving 
the problems (P) and (D) by fully exploiting the dual block-angular structure.
Our method is based on the abstract inexact symmetric Gauss-Seidel proximal ADMM (sGS-ADMM) framework  developed in \cite{CST}. Specifically, we develop 
a concrete sGS-ADMM to solve the block-angular problem (P) via its dual problem (D). 
While our algorithm is based on the abstract algorithmic framework in \cite{CST}, 
the key issue in the current algorithmic design is 
to fully exploit the block-angular structure in our target problem  to 
decompose the computation of the variable blocks, with the goal of making
the subproblems in each iteration easier to solve. 
Moreover, we also need to develop the corresponding numerical schemes to solve the
subproblems efficiently. 
As we shall see later, our method can fully exploit
the block-angular structure in the problem and hence the subproblems in each iteration 
can be solved efficiently. In addition, it also offers us the flexibility 
to solve the subproblems in parallel.

Now we summarize our contributions in this paper.
\begin{enumerate}
	\item 
	We derived concrete DBA convex composite conic programming problems arising from
	DNN relaxations of common mixed integer programming problems such as
	uncapacitated facility location problems.
	
	\item {We proposed and implemented a specialized algorithm to solve large scale DBA problems {via their duals}. In particular, by leveraging on the recently developed sGS decomposition theorem \cite{LST-sGS}, we have designed proximal ADMM and proximal ALM frameworks that are able to fully exploit the underlying block angular structure. 
		This is done by designing appropriate proximal terms within 
		a proximal ADMM or ALM framework to enable the decomposition of 
		the dual block variables' computation. By further designing novel numerical schemes to efficiently solve
		the decomposed subproblems, each iteration of the overall algorithm 
		is relatively easy and inexpensive to execute.  
		{In addition, 
			our algorithm can conveniently take advantage of parallel computation to further improve its
			computational efficiency.}
	}

	\item We conduct comprehensive numerical experiments to compare our algorithm with state-of-the-art solvers. The results show that our algorithm outperforms the others especially when the constraint set $\cK_i$ is the cone $\S^{n_i}_+$ or the quadratic coefficient matrix in the objective function is dense. For two stage stochastic linear programming, we also compare our algorithm with the progressive hedging algorithm (PHA) in \cite{RTRW}. The numerical results show that our algorithm can perform much better when the number of scenarios is large.
\end{enumerate}

In the current literature, {as far as we are aware of, all the popular} augmented Lagrangian based decomposition methods (including the PHA) for DBA linear programming
problems directly target the primal problem (P). Here our methods are designed
to deal with the structure in the dual problem (D) while at the same time solving the primal
problem (P) as a by-product. The key advantage in this dual approach is that we do not need to introduce
any extra variables and constraints, unlike the case of directly dealing with the primal
problem where a large number of extra variables and constraints must be introduced to
split the first-stage variable $x$ so as to decompose the computation of the
second-stage variables $\bar{x}_i$ ($i = 1,\ldots,N$). 
In contrast, for our dual approach, the decomposition
of the computation of the dual variables can be obtained naturally with the sGS
decomposition framework.

This paper is organized as follows. We will describe some example classes of {convex} programming problems with dual block-angular structures in section 2. In section 3, we will derive the dual (D) of the primal problem (P). We provide some preliminaries in section 4 before 
providing the algorithmic details of our inexact symmetric Gauss-Seidel proximal ADMM and ALM for solving (D) in section 5. In section 6, we discuss an extension where we incorporate a semismooth Newton-CG method in our proposed algorithm. We conduct numerical experiments in section 7 to evaluate the performance of our algorithms. Finally, we conclude the paper in the last section.

\bigskip 
\textbf{Notation.}
\begin{itemize}
	\item \blue{$\mathbb R^d$ denotes the $d$-dimensional Euclidean space, $\mathbb S^d$ denotes the space of the $d$-dimensional symmetric matrices, $\mathbb S^d_+$ denotes the $d$-dimensional positive semidefinite cone.}
	\item {Let $\cT : \cX\rightarrow\cX$ be a self-adjoint positive semidefinite map. We denote $\|x\|_\cT := \sqrt{\inprod{x}{\cT x}}$ as the semi-norm induced by $\cT$.
	For any linear map ${\cal B} :\cX \to \cY$, where $\cX$ and $\cY$ are finite-dimensional 
	inner product spaces, its induced operator norm $\norm{\cB} = \max_{x\not=0} \norm{\cB x}/\norm{x}.$
	}
	\item Let $f :\cX\rightarrow(-\infty,+\infty]$ be a closed proper convex function. Then, we denote dom$f$ and $\partial f$ as its effective domain and subdifferential, respectively. 
	\item We denote $f^*$ {as} the Fenchel conjugate function of a convex closed function $f$, which is defined {by} $f^*(x):=\sup_y\{\inprod{x}{y}-f(y)\}$.
	\item The Moreau-Yosida envelope of $f$ is defined by $M_f(x):=\min_y\{f(y)+\frac{1}{2}\|x-y\|^2\}$, and proximal mapping of $f$ is defined by $\text{Prox}_f(x):=\arg\min_y\{f(y)+\frac{1}{2}\|x-y\|^2\}$.
	\item The following formulas are frequently used in this \blue{paper}:
	\[
	\begin{aligned}
	&\grad M_{f^*/t}(x)=t^{-1}\text{Prox}_{tf}(tx)\quad and\quad \grad M_{tf}(x)=t\text{Prox}_{f^*/t}(x/t)\\
	&x=\text{Prox}_{tf}(x)+t\text{Prox}_{f^*/t}(x/t)\quad and\quad \frac{1}{2}\|x\|^2=M_{tf}(x)+t^2M_{f^*/t}(x/t)
	\end{aligned}
	\] 
	\item For given square matrices $P_1,\ldots,P_\ell$, we use ${\rm diag}(P_1,\ldots,P_\ell)$ to denote the block diagonal
	matrix with diagonal blocks given by $P_1,\ldots,P_\ell$. 
\end{itemize}

\subsection{Related literature}

{
	DBA linear-quadratic problems are well-known in the context of two-stage stochastic programming. In the pioneering work of \cite{RTRW}, a multistage programming problem was considered where the variables are vectors in some feasible set in the Euclidean spaces, and the objective function is linear-quadratic. In that paper, the progressive hedging algorithm (PHA) 
	was proposed, which is an algorithmic framework based on
	{a proximal augmented Lagrangian method whereby at each iteration, 
		each scenario requires the computation of one linear-quadratic programming subproblem. The PHA can be costly
		when the number of scenarios is large and/or when the subproblems are difficult to solve.}
	Subsequently, augmented Lagrangian based methods like Diagonal Quadratic Approximation method (DQA) in \cite{ruszcz95} and Accelerated Distributed Augmented Lagrangian method (ADAL) in \cite{CDZ} have been proposed. These two methods are initially designed to solve large-scale convex optimization problems with primal block angular (PBA) structure that have many linking constraints. Both methods aim to make their subproblems easy to solve by modifying the augmented Lagrangian function and try to make it separable. However, most of these methods are still difficult to be applied to a large-scale problem when the underlying feasible set is nontrivial. {Just like the PHA, the computational cost of solving those separable smaller linear-quadratic programming subproblems can become expensive when the number of subproblems is 
		large and/or when the subproblems are difficult to solve.}
	
	Other than two-stage stochastic programming problem, there are less literature devoted to designing specialized methods for directly solving dual block angular problems. Most of them construct a general framework and solve DBA problems as a by-product. For example, in \cite{Kontogiorgis}, the authors developed a framework based on  alternating directions method for solving the minimization of the sum of two convex functions subject to linear constraints. It is then applied successfully to solve the Fermat-Weber problem recasted into a dual-block angular problem. 
}

While many papers have been written on augmented
Lagrangian based decomposition methods for solving DBA linear and quadratic programming
problems, there are far fewer papers on other classes of DBA problems such as DBA convex
conic programming problems where the cones are more general than the nonnegative orthants. In
\cite{SM1}, {the authors} proposed a log-barrier decomposition interior-point method (IPM) for a 2-stage stochastic
SDP that is expressed as the dual of the primal block-angular form, wherein the algorithm
is based on the pioneering work of Zhao \cite{zhao2001} on a log-barrier Benders decomposition IPM for
a 2-stage stochastic LP. But note that such a decomposition method is not directly designed for
the DBA problem (1), so we shall not touch on such a method any further.

In \cite{SM2}, a log-barrier Benders decomposition IPM is proposed for a 2-stage stochastic convex
QP expressed in the dual block-angular form, and the polynomial iteration complexity of the
method is analysed, {but there are no numerical implementation.} The IPM is extended in \cite{Chen-Mehrotra} for a 2-stage stochastic convex optimization
problem where the nonnegative orthants are replaced by compact convex sets whose barrier
functions (and the corresponding gradients and Hessians) are assumed to be readily computable. 

Our algorithms are applicable to  the very general standard DBA problem (1). In \cite{Chen-Mehrotra}, a
theoretical log-barrier decomposition IPM has been proposed to solve dual block-angular
convex programming problems with linear objective functions (i.e, with $\theta(\cdot)\equiv 0$ and $\bar{\theta}(\cdot)\equiv 0$). However, the implementation and practical evaluation of the methods so far
have not been performed.

\blue{
In \cite{LST-PBA}, the authors considered convex optimization problems with the
primal block angular (PBA) structure, which has a different structure from the dual block-angular structure considered in this paper. They also consider the dual problem to deal with the PBA structure by sGS-ADMM, but the design of the proximal terms are different from ours in this paper. Moreover, numerical experiments on more general cone like SDP cone are not performed in \cite{LST-PBA}.}

\section{Example classes of \eqref{eq-DBA}}

\subsection{Examples from 2-stage stochastic {convex} programming problems}

Here we show that  the dual block-angular structure shown in 
\eqref{eq-DBA} naturally arises from a 2-stage stochastic {convex} programming problem.
Consider the following 2-stage stochastic optimization problem:
\begin{eqnarray}
\min_{x\in\cX} \Big\{ \theta(x) + \inprod{c}{x} + \blue{E_\xi \left(Q(x;\xi)\right)}  \mid
A x \;=\; b, \; x \in \cK \Big\}.
\label{eq-stoc}
\end{eqnarray}
Here $\xi$ is a random vector, $E_\xi(\cdot)$ denotes the expectation with respect to $\xi$, and
\begin{eqnarray*}
	Q(x;\xi) \;=\; \min_{\tx\in\cX_1} \left\{\tilde{\theta}_\xi(\tx)+ \inprod{\tc_\xi}{\tx} \,:\,
	\tB_\xi \tx= \tb_\xi - B_\xi x, \; \tx \in \cK_1 \right\},
\end{eqnarray*}
where for each given scenario $\xi$, 
$\tc_\xi\in \cX_1$, $\tb_\xi\in \cY_1$ are given data, 
$B_\xi:\cX\to \cY_1$,
$\tB_\xi :\cX_1\to \cY_1$ are given linear maps, and $\cK_1 \subset \cX_1$ is 
a given closed convex set, 
$\tilde{\theta}_\xi:\cX_1 \to (-\infty,\infty]$ is a given proper closed convex function. 
Note that there is no loss of generality in considering only equality constraints
in \eqref{eq-stoc} since constraints of the form 
$b- A x \in \cQ$, $x\in \cK$ can always be reformulated as 
$[A,I](x;s) = b$, $(x;s) \in \cK\times \cQ$. 

By sampling $N$ scenarios for $\xi$, one may 
approximate $E_\xi Q(x;\xi)$ by the sample mean $\sum_{i=1}^N p_i Q(x;\xi_i)$, 
with $p_i$ being the probability of occurrence of the $i$th scenario,
and hence
approximately solve
\eqref{eq-stoc} via the following deterministic optimization problem
which has exactly the same form as \eqref{eq-DBA}:
\begin{scriptsize}
\begin{eqnarray}
\min \left\{ \theta(x) + \inprod{c}{x} + \sum_{i=1}^N \blue{\left(\btheta_i(\bx_i)+\inprod{\bc_i}{\bx_i}\right)}  \Bigm\vert \begin{array}{l} Ax = b, \; x\in \cK\\[5pt] B_i x + \bB_i\by_i = \bb_i, \;\by_i\in\cK_i,\; \forall\; i=1,...,N
\end{array} \right\}\;\label{eq-stoc2}
\end{eqnarray}
\end{scriptsize}
where $\bc_i=p_i \tc_{\xi_i}$ and $\btheta_i = p_i \tilde{\theta}_{\xi_i}$,
$B_i = B_{\xi_i}, \bB_i=\tB_{\xi_i}, \bb_i=\tb_{\xi_i}$ are the data,
and $\bx_i=\tx_{\xi_i}\in\cX_1$ is the
second stage decision variable associated with the $i$th scenario.

\subsection{Examples from doubly nonnegative relaxations of uncapacitated facility location problems}

The uncapacitated facility location (UFL) problem is one of the most studied 
problems in operations research. Here we consider the general 
UFL with linear and/or separable convex quadratic allocation costs
that was introduced in \cite{Gunluk}:
\begin{eqnarray}
\begin{array}{rl}
\mbox{(UFL)} \quad 
\min & \sum_{i=1}^p c_i u_i + \sum_{i=1}^p \sum_{j=1}^q \blue{\left(p_{ij} s_{ij} + \frac{1}{2}q_{ij} s_{ij}^2\right)}
\\[5pt]
\mbox{s.t} & \sum_{i=1}^n s_{ij} = 1, \;\forall\;  j=1,\ldots,q
\\[5pt]
&  0 \leq s_{ij} \leq u_i,  \; \forall\; i=1,\ldots,p, \; j=1,\ldots,q
\\[5pt]
& u_i \in \{0,1\} \; \forall\; i=1,\ldots,p,
\end{array}
\label{eq-UFL}
\end{eqnarray}
where $c_i\geq 0$ for all $i$, and $\blue{p_{ij}},q_{ij} \geq 0$ for all $i,j$ are given data. 
Observe that the allocation cost for customer $j$ is 
a convex quadratic function of his demand $s_{ij}$ being served by the
opened facility $i$. 

Let $U = uu^T$. One can see that the constraint $u_i \in \{0,1\}$ 
is equivalent to the constraint that  $u_{i}^2 = U_{ii} = u_{i}$.
Also, by introducing the nonnegative slack variable $z_{ij}$, we can convert the inequality 
constraint $s_{ij}\leq u_i$ to the equality constraint $s_{ij} + z_{ij} = u_i$. 

Let $S = (s_{ij})\in \RR^{p\times q}$ and $Z = (z_{ij})\in \RR^{p\times q}$. 
Also let $e\in\RR^p$ be the vector of all ones 
and $S_j$ be the $j$th column of $S$.  
Then one can express the problem \eqref{eq-UFL} equivalently as follows:
\begin{eqnarray*}
	\begin{array}{rl}
		\min & \sum_{i=1}^p c_i u_i + \sum_{i=1}^p \sum_{j=1}^q \blue{\left(p_{ij} s_{ij} + \frac{1}{2}q_{ij} s_{ij}^2\right)}
		\\[5pt]
		\mbox{s.t} & e^T S_j = 1, \;\;\forall\;  j=1,\ldots,q
		\\[5pt]
		&  S_j + Z_j = u, \;\; S_j, Z_j \geq 0,\;\; \forall\; j=1,\ldots,q
		\\[5pt] 
		& u- \diag{U} = 0,   \quad U = uu^T,  \; u\geq 0, \; U\in\S^p.
	\end{array}
\end{eqnarray*}
Let $\blue{P = (p_{ij})}\in \RR^{p\times q}$ and $Q=(q_{ij})\in\RR^{p\times q}.$ 
By relaxing the rank-1 constraint $U = uu^T$ to $U\succeq uu^T$ and using the
equivalence that $U\succeq uu^T$ if and only if $[1,\; u^T; u ,\; U]\succeq 0$, we get the following 
doubly nonnegative (DNN) relaxation of  \eqref{eq-UFL}: 
\begin{eqnarray}
\begin{array}{rl}
\min & \inprod{[0; c]}{[\alp; u]} + \sum_{j=1}^q \blue{\left(\inprod{P_j}{S_j} + \frac{1}{2} \inprod{S_j}{\diag{Q_j}S_j}\right)}
\\[8pt]
\mbox{s.t} &
\left[\begin{array}{c} 0 \\ u \end{array}\right] 
+
\left[\begin{array}{cc}
e^T &  0^T \\
-I_p         & -I_p
\end{array}\right] \left[\begin{array}{c} S_j \\ Z_j \end{array}\right]
= 
\left[\begin{array}{c} 1 \\ 0 \end{array}\right], \quad j = 1,\ldots q
\\[15pt]
& u - \diag{U} = 0, \;\;  \alp = 1,
\\[5pt]
& \mbU := \left[\begin{array}{cc} \alp & u^T \\[3pt] u & U \end{array}\right] \in \S^{1+p}_+, \;\;
\mbU \geq 0,\;\; S_j, Z_j \geq 0, \; j=1,\ldots, q.
\end{array}
\label{eq-UFL-DNN}
\end{eqnarray}
Now we can express the above DNN  programming problem  in the form \eqref{eq-DBA} by
defining the following inner product spaces, functions, cones and linear maps:
\begin{eqnarray*}
	& \cX = \S^{1+p},  \;\; \cY = \RR^{1+p}, \;\; \cX_j = \RR^{2q},\;\; \cY_j = \RR^{1+p},
	\; j=1,\ldots,q,  &
	\\[5pt]
	&\theta(\mbU) = \delta_{\cN}(\mbU), \; \cK = \S^{1+p}_+ ,& \\[5pt]
	&\btheta_j([S_j; Z_j]) = \frac{1}{2}\inprod{S_j}{\diag{Q_j}S_j}, 
	\; \cK_j = \RR^{2q}_+,\; j=1,\ldots,q, &
	\\[5pt]
	& b= \left[\begin{array}{c} 1 \\[5pt] 0 \end{array}\right] \in \RR^{1+p}, \;\;
	A(\mbU) = \left[\begin{array}{c} \alp \\[5pt] u-\diag{U}\end{array}\right]
	=\left[\begin{array}{c}
		\inprod{A_0}{\mbU} \\
		\inprod{A_1}{\mbU}\\
		\vdots \\
		\inprod{A_p}{\mbU}
	\end{array} \right],
	&
	\\[5pt]
	& 
	B_j(\mbU) = \left[ \begin{array}{c} 0 \\[5pt] u \end{array}\right]
	=
	\left[\begin{array}{c}
		\inprod{O}{\mbU} \\
		\inprod{E_1}{\mbU}\\
		\vdots \\
		\inprod{E_p}{\mbU}
	\end{array} \right],
	\;\;
	\bB_j = \left[\begin{array}{cc}
		e^T &  0^T \\
		-I_p         & -I_p
	\end{array}\right] \in  \RR^{(1+p)\times2q} ,& \\[5pt]
	&\bb_j= \left[\begin{array}{c} 1 \\ 0 \end{array}\right] \in \RR^{1+p}, \;\; j=1,\ldots, q,
	&
\end{eqnarray*}
where $\cN = \{ X \in \S^{1+p} \mid X \geq 0\}$, $Q_j$ is the $j$th column of $Q$, and
\begin{eqnarray*}
	A_0 = \left[\begin{array}{cc} 
		1  & 0^T \\[3pt]  0  & 0_{p\times p}
	\end{array} \right],\quad 
	A_i = \left[\begin{array}{cc} 
		0  & \frac{1}{2} e_i^T \\[3pt] \frac{1}{2} e_i  & -e_i e_i^T 
	\end{array} \right],\quad 
	E_i =  \left[\begin{array}{cc} 
		0  & \frac{1}{2} e_i^T \\[3pt] \frac{1}{2} e_i  & 0_{p\times p}
	\end{array} \right], \;\; i=1,\ldots,p. 
\end{eqnarray*} 
In the above, $e_i$ is the $i$th unit vector in $\RR^p.$
Observe that for the problem \eqref{eq-UFL-DNN}, we have that
$$
B_1  = \cdots = B_q, \quad \bB_1 = \cdots = \bB_q.
$$
Thus given any $\mbU\in \cX$ and $\by\in \cX_1\times\cdots\times \cX_q$, the evaluation of $B(\mbU)$ and $B^*\by$ can be 
done very efficiently since $B(\mbU) = [B_1(\mbU)]_{j=1}^q$ and
$B^*\by = B_1^*(\by_1 + \cdots + \by_q)$.
For later purpose, we note that $\bB_j \bB_j^*$ has a very simple inverse given by
\begin{eqnarray}
(\bB_j \bB_j^*)^{-1} &=& \left[\begin{array}{cc} 
0 & 0^T \\ 0 &\frac{1}{2}I_p 
\end{array}\right] +
\frac{1}{2p} \left[\begin{array}{c} 
2 \\ e 
\end{array}\right] [2,\; e^T], \quad j=1,\ldots, q.
\label{eq-UFL-bBbBt}
\end{eqnarray}

\begin{remark}
If the problem \eqref{eq-UFL} has additional constraints in $u$, say $Hu = \beta$, then the 
corresponding constraints in the 
DNN relaxation
of the problem would become
\begin{eqnarray*}
\inprod{H_k} {\mbU} = \beta_k, \quad \inprod{\blue{\widehat{H}_k^{}}}{\mbU} = \beta_k^2,\quad
\quad k=1,\ldots m
\end{eqnarray*}
where given the $k$th row $h_k^T$ of $H$,
$$
H_k =  \left[\begin{array}{cc} 
0 & h_k^T/2 \\[3pt] h_k/2 & 0
\end{array}\right], \quad 
\widehat{H}_k =  \left[\begin{array}{cc} 
0 & 0^T \\[3pt] 0 & h_k h_k^T
\end{array}\right], \quad k=1,\ldots,m.
$$
\end{remark}
\section{The derivation of (D) and KKT conditions}
\label{sec:dual}

By introducing the auxiliary variables, $u=x$, $\bar{u}=\bx$,
$s=x$, $\bar{s}=\bx$
and replacing 
$\theta(x)$, $\btheta(\bx)$, 
$\delta_{\cK}(x)$, $\delta_{\bcK}(\bx)$
by $\theta(u)$, $\btheta(\bar{u})$,
$\delta_{\cK}(s)$, $\delta_{\bcK}(\bar{s})$
respectively, 
one can derive the dual of (P). 
Consider the Lagrangian function for (P):
\begin{eqnarray*}
	&& \hspace{-0.7cm}
	\cL(x,\bx,u,\bar{u},s,\bar{s};y,\by,v,\bv,z,\bz) \\
	&:=& \theta(u) + \inprod{c}{x} 
	+ \btheta(\bar{u})+\inprod{\bc}{\bx} + \delta_{\cK}(s) + \delta_{\bcK}(\bar{s}) - \inprod{y}{Ax-b} - \inprod{\by}{Bx+\bB \bx - \bb} \\
	&& - \inprod{v}{x-u} - \inprod{\bv}{\bx - \bar{u}} - \inprod{z}{x-s} - \inprod{\bz}{\bx - \bar{s}}\\
	&=& \inprod{c-A^* y - B^* \by - v - z}{x} + \blue{\inprod{\bar c-\bB^* \by - \bv - \bz}{\bx}} + \theta(u) + \inprod{v}{u} + \btheta(\bar{u}) + \inprod{\bv}{\bar{u}} \\
	&& + \delta_{\cK}(s) + \inprod{z}{s} + \delta_{\bcK}(\bar{s}) + \inprod{\bz}{\bar{s}} + \inprod{b}{y} + \inprod{\bb}{\by},
\end{eqnarray*}
where $x,u,s,v,z \in \cX$, $\bx,\bar{u},\bar{s},\bv,\bz \in\bcX$, $y\in \cY$ and $\by \in \bcY$. We have the following:
\blue{
\begin{eqnarray*}
	\inf_x \cL(x,\bx,u,\bar{u},s,\bar{s};y,\by,v,\bv,z,\bz) &\Longleftrightarrow& 
	\inf_x \{\inprod{c-A^* y - B^* \by - v - z}{x}\} \\
	&=& \begin{cases}
		0, & \text{if }  c-A^* y - B^* \by - v - z = 0, \\
		-\infty, & \text{otherwise}.
	\end{cases} \\
	\inf_{\bx} \cL(x,\bx,u,\bar{u},s,\bar{s};y,\by,v,\bv,z,\bz) &\Longleftrightarrow& 
	\inf_{\bx}\{\inprod{\bc-\bB^* \by - \bv - \bz}{\bx}\} \\
	&=& \begin{cases}
		0, & \text{if }  \bB^* \by + \bv + \bz = \bc, \\
		-\infty, & \text{otherwise}.
	\end{cases} \\
	\inf_u \cL(x,\bx,u,\bar{u},s,\bar{s};y,\by,v,\bv,z,\bz) &\Longleftrightarrow& 
	\inf_u \{\theta(u) + \inprod{v}{u}\} = -\theta^*(-v) \\
	\inf_{\bar{u}}\cL(x,\bx,u,\bar{u},s,\bar{s};y,\by,v,\bv,z,\bz) &\Longleftrightarrow&
	\inf_{\bar{u}}\{\btheta(\bar{u}) + \inprod{\bv}{\bar{u}}\} = -\btheta^*(-\bv)\\
	\inf_s \cL(x,\bx,u,\bar{u},s,\bar{s};y,\by,v,\bv,z,\bz) &\Longleftrightarrow& 
	\inf_u \{\delta_{\cK}(s) + \inprod{z}{s}\} = - \delta_{\cK}^* (-z)\\
	\inf_{\bar{s}}\cL(x,\bx,u,\bar{u},s,\bar{s};y,\by,v,\bv,z,\bz) &\Longleftrightarrow& 
	\inf_{\bar{s}}\{\delta_{\bcK}(\bar{s}) + \inprod{\bz}{\bar{s}}\} = -\delta_{\bcK}^*(-\bz),
\end{eqnarray*}}
\blue{where "$\Longleftrightarrow$" means equivalence by removing the irrelevant constants.}

Hence the dual problem of (P) is given by
\begin{eqnarray*}
	&\mbox{(D)}& \max_{y,\by,v,\bv,z,\bz} \; \inf_{x,\bx,u,\bar{u},s,\bar{s}}\cL(x,\bx,u,\bar{u},s,\bar{s};y,\by,v,\bv,z,\bz) 
	\\
	&&
	\begin{array}{rl} 
		=\;\max_{y,\by,v,\bv,z,\bz} & -\theta^*(-v)-\btheta^*(-\bv)- \delta_{\cK}^* (-z)-\delta_{\bcK}^*(-\bz) + \inprod{b}{y} + \inprod{\bb}{\by} 
		\\[5pt]
		\mbox{s.t.} &
		c-A^* y - B^* \by - v - z = 0, 
		\\[5pt]
		& \bar c - \bB^* \by - \bv - \bz = 0.
	\end{array} 
\end{eqnarray*}
\blue{
Consider the Lagrangian function for (D),
\[
\begin{aligned}
&\cL(v,\bar v,z,\bar z,y,\bar y;x,\bar x)\\:=&
\;\;
\theta^*(-v)+\bar\theta^*(-\bar v)+\delta^*_{\cK}(-z)+\delta_{\bar\cK^*}(-\bar z)-\inprod{b}{y}-\inprod{\bar b}{\bar y}\\&+\inprod{x}{c-A^*y-B^*\bar y-v-z}+\inprod{\bar x}{\bar c-\bar B^*\bar y-\bar v-\bar z}.
\end{aligned}
\]}
Assume that the feasible regions of both the primal and dual problems have nonempty interiors. Then the optimal solutions for both problems exist, and they satisfy
the following Karush-Kuhn-Tucker (KKT) conditions for the problems (P) and (D):
\blue{
\begin{eqnarray}
\left\{
\begin{array}{ccc}
& Ax = b, \quad Bx + \bB\bx = \bb, & \\[3pt]
& A^*y + B^* \by + z + v =c,\quad \bB^*\by + \bz + \bar{v} = \bc, & \\[5pt]
& 0\in \partial\delta_\cK^*(-z)-x \Leftrightarrow -z={\rm Prox}_{\delta_\cK^*}(x-z)\Leftrightarrow x = \Pi_{\cK}(x-z), \\ 
&0\in \partial\delta_{\bar\cK^*}(-\bar z)-\bar x \Leftrightarrow -\bar z={\rm Prox}_{\delta_{\bar\cK^*}}(\bar x-\bar z)\Leftrightarrow \bar x = \Pi_{\bar\cK}(\bar x-\bar z), \\
&0\in \partial\theta^*(-v)-x \Leftrightarrow -v={\rm Prox}_{\theta^*}(x-v)\Leftrightarrow x = {\rm Prox}_{\theta}(x-v), \\
&0\in \partial\bar\theta^*(-\bar v)-\bx \Leftrightarrow -\bar v={\rm Prox}_{\bar\theta^*}(\bx-\bar v)\Leftrightarrow \bar x = {\rm Prox}_{\bar\theta}(\bx-\bar v).\\
\label{eq-KKT}
\end{array}
\right.
\end{eqnarray}
}
\blue{where $\Pi_{\cK}(\cdot)$ and $\Pi_{\bcK}(\cdot)$ denote the projection mappings onto
$\cK$ and $\bcK$, respectively. }

\section{Preliminaries}

{In this section, we will first present the sGS decomposition theorem, which is important for formulating the main algorithmic framework that is given in \cite{CST}. In short, this theorem helps to decompose a multi-block problem into smaller subproblems with additional inexpensive steps but not sacrificing the convergence guarantee. Then, we will present the abstract algorithmic template by fitting our problem in the context of \cite{CST}. }

\subsection{sGS decomposition theorem}
Let $\cX_1,\ldots,\cX_s$ ($s\geq 2$) be given finite-dimensional inner product spaces,  
$\cQ : \cX\to\cX$ be a given self-adjoint positive semidefinite linear operator, and
$c\in \cX$, where $\cX = \cX_1\times\cdots\times \cX_s.$
Suppose that $\cQ$ is partitioned according to $\cX_1\times\cdots\times \cX_s$ as
\begin{eqnarray}
\cQ &=& \left[ \begin{array}{ccc}
Q_{1,1} &\dots  & Q_{1,s}  \\[5pt]
\vdots           &\vdots         & \vdots                \\[5pt]
Q_{1,s}^* &\dots  & Q_{s,s}
\end{array}\right]  \;=\; \cU + \cD+\cU^*,
\label{eq-Q}
\end{eqnarray}
where $\cD$ and $\cU$ denote the block diagonal part and the strictly upper block triangular part of 
$\cQ$, respectively. That is,
\begin{eqnarray}
\label{eq-UD}
\cU =
\left[ \begin{array}{cccc}
\mb0 & Q_{1,2} &\dots & Q_{1,s} \\
& \ddots   &        & \vdots \\
&              & \ddots & Q_{s-1,s} \\[5pt]
&              &            &\mb0
\end{array}\right],
\quad \cD =  \left[ \begin{array}{cccc}
Q_{1,1} &  & \\[3pt]
& Q_{2,2}   &        & \\
&              & \ddots & \\
&              &            &Q_{s,s}
\end{array}\right].
\end{eqnarray}
Assuming that $Q_{ii}$ is positive definite for all $i=1,\ldots,s$, then we can define the 
following symmetric Gauss-Seidel (sGS) linear operator associated with $\cQ$: 
\begin{eqnarray}
{\bf sGS}(\cQ) = \cU \cD^{-1} \cU^*. 
\label{eq-sGS}
\end{eqnarray}
For a given $x= (x_1;\ldots; x_s)$, we define
\begin{eqnarray*}
	x _{\geq i} = (x_i;\ldots; x_s), \quad x_{\leq i} = (x_1;\ldots;x_i), \quad i = 1,\ldots,s.
\end{eqnarray*}
We also define $x_{\geq s+1} = \emptyset$.
We have the following theorem which will be useful for us to design
a decomposition method to solve (D).

\begin{theorem} \mbox{\rm ({\bf sGS decomposition theorem \cite[Theorem 1]{LST-sGS}})} Let $ z\in\cX$ be given. {Assume that $\cQ\succeq0$ and $Q_{i,i}$ are positive definite.} Suppose
	$\delta', \delta\in\cX$ are two given error
	vectors with $\delta_1' = \delta_1$. Define
	\begin{eqnarray}\label{eq-Delta}
	\Delta(\delta',\delta):= \delta + \cU \cD^{-1}(\delta -\delta').
	\end{eqnarray}
	Let $p$ be a closed proper convex function, and $q(x) =  \frac{1}{2}\inprod{x}{\cQ x} 
	- \inprod{c}{x}\;\forall\; x\in \cX.$
	Consider the following  convex composite quadratic programming problem:
	\begin{eqnarray}
	x^+ := \argmin \;\Big\{ p(x_1) +q(x) + \frac{1}{2} \norm{x-z}_\cS^2 
	-\inprod{x}{\Delta(\delta',\delta)}
	\;\mid\; x\in \cX
	\Big\}  ,
	\label{prox-QP}
	\end{eqnarray}
	where $\cS = {\bf sGS}(\cQ)$ is the sGS linear operator 
	associated with $\cQ$. 
	Then $x^+$ can be computed in the following symmetric Gauss-Seidel fashion. 
	For $i=s, \ldots,2,$  compute $x_i^\prime\in \cX_i$ defined by
	\begin{equation}
	\begin{aligned}
	x^\prime_i :={}& \argmin_{x_i\in \cX_i} \; \big\{ p(z_1) + q(z_{\le i-1};x_i;x^\prime_{\ge i+1})-\inprod{\delta_i'}{x_i} \big\}
	\\[5pt]
	={}& Q_{i,i}^{-1}\big(c_i +\delta_i' - \mbox{$\sum_{j=1}^{i-1}$} Q_{j,i}^*z_j - \mbox{$\sum_{j=i+1}^s$} Q_{i,j}x^\prime_j \big).
	\end{aligned}
	\label{xpimrei}
	\end{equation}
	Then the optimal solution $x^+$ for \eqref{prox-QP} can be computed exactly via the following steps:
	\begin{equation} \label{prox-nT}
	\left\{
	\begin{aligned}
	x_1^+ ={}& \argmin_{x_1\in\cX_1}\; \big\{p(x_1) + q(x_1;x^\prime_{\ge 2}) - \inprod{\delta_1}{x_1} \big\}
	\\[5pt]
	x_i^+ ={}& \argmin_{x_i\in\cX_i} \;\big\{ p(x_1^+) + q(x^+_{\le i-1};x_i;x^\prime_{\ge i+1}) - \inprod{\delta_i}{x_i} \big\}\\[5pt]
	={}& Q_{i,i}^{-1}\big(c_i+\delta_i -\mbox{$\sum_{j=1}^{i-1}$}
	Q_{j,i}^*x_j^+ - \mbox{$\sum_{j=i+1}^s$} Q_{i,j}x^\prime_j\big),\quad i=2,\ldots,s.
	\end{aligned}
	\right.
	\end{equation}
	\label{sGS-decom}
\end{theorem} 

Note that the role of the error vectors
$ \delta'$ and $ \delta$ in the above block  sGS decomposition theorem is to 
allow for inexact computation in \eqref{xpimrei} and \eqref{prox-nT}. 
There is no need to know these error vectors in advance. We should  view 
the computed $x_i'$ and $x_i^+$
from \eqref{xpimrei} and \eqref{prox-nT}
as  approximate solutions to the minimization subproblems 
without the terms involving $\delta_i'$ and $\delta_i$. Once these approximate solutions have been computed, the errors incurred would generate the vectors $ \delta_i'$ and $ \delta_i$ automatically.   With these known error vectors, we know that the computed approximate solutions are 
the exact solutions to the slightly perturbed 
minimization problems in \eqref{xpimrei} and \eqref{prox-nT}.
In particular, for $i=s,\ldots,2$, we know that $\delta_i'$  is the 
residual vector obtained when we solve the linear system of equations $Q_{ii} x^\prime_i = 
c_i - \mbox{$\sum_{j=1}^{i-1}$} Q_{j,i}^*z_j - \mbox{$\sum_{j=i+1}^s$} Q_{i,j}x^\prime_j$
in \eqref{xpimrei} approximately.
{Similarly,} $\delta_i$ ($i=2,\ldots,s$) is the residual vector obtained when we solve the 
linear system of equations in \eqref{prox-nT}.

\subsection{An abstract inexact proximal ADMM template}

To make use of the abstract inexact proximal ADMM
template provided in \cite{CST}, we first write (D) in the form that was presented in that paper:
\begin{eqnarray}\label{eq-2blksCCP}
\min_{\mu\in\ctX,\nu\in\ctY} \big\{p(\mu_1) + f(\mu_1,\mu_2) + q(\nu_1) + g(\nu_1,\nu_2) \; \mid \; \cF^* \mu + \cG^* \nu = h \big\},
\end{eqnarray}
where $\ctX = (\bcX\times \cX)\times \cY$, $\ctY = (\cX\times \bcX)\times\bcY$,
\begin{eqnarray*}
	&&\mu_1:=[\bz;z], \quad \mu_2 = y, \quad \nu_1 = [v;\bv], \quad \nu_2 = \by, 
	\\[5pt]
	&&p(\mu_1) := \delta_{\bcK}^*(-\bz) + \delta_{\cK}^* (-z), \quad 
	f(\mu_1,\mu_2) := - \inprod{b}{y}, \\[5pt]
	&&q(\nu_1) := \theta^*(-v) + \btheta^*(-\bv), \quad 
	g(\nu_1,\nu_2) := -\inprod{\bb}{\by}, \\[5pt]
	&&\cF^* := \begin{bmatrix}
		0 & I & \cA^* \\ I & 0 & 0
	\end{bmatrix}, \quad 
	\cG^* := \begin{bmatrix}
		I & 0 & B^* \\ 0 & I & \bar{B}^*
	\end{bmatrix}, \quad 
	h := \blue{\begin{bmatrix}
		c \\ \bc
	\end{bmatrix}}.
\end{eqnarray*}

Given a penalty parameter $\sig > 0$, 
the augmented Lagrangian function associated with \eqref{eq-2blksCCP} is given as follows. For $\mu :=(\mu_1,\mu_2)\in (\bcX\times \cX)\times \cY$, 
$\nu := (\nu_1,\nu_2) \in (\cX\times \bcX)\times\bcY$, $\omega = [x;\bx]\in\cX\times \bcX$,
\begin{eqnarray*}
	\cL_\sig(\mu,\nu;\omega) 
	&=& p(\mu_1) + f(\mu_1,\mu_2) + q(\nu_1) + g(\nu_1,\nu_2) -\frac{1}{2\sig}\norm{\omega}^2 \\
	&& + \frac{\sig}{2}\norm{\cF^*\mu + \cG^*\nu-h+\sig^{-1}\omega}^2.
\end{eqnarray*}
\blue{Let $\partial_\mu\cL_\sigma(\mu,\nu;\omega)$ and $\partial_\nu\cL_\sigma(\mu,\nu;\omega)$ denote the $\mu$ and $\nu$ components of the subdifferential respectively.}
The inexact proximal ADMM template
for solving \eqref{eq-2blksCCP} is presented next.

\fbox{\parbox{\textwidth}{
		\textbf{Algorithm imsPADMM in \cite{CST}:} An abstract inexact proximal ADMM for solving \eqref{eq-2blksCCP}. \\
		Let $\tau \in (0,(1+\sqrt{5})/2)$ be the step-length and $\{\tilde{\eps}_k\}_{k\ge 0}$ be a summable sequence of nonnegative numbers. Choose an initial point $(\mu^0,\nu^0;\omega^0)\in ({\rm dom}p\times \cY) \times ({\rm dom}q \times \bcY) \times (\cX\times \bcX)$. Let $\tilde{\cS},\tilde{\cT}$ be positive semidefinite linear operators that are to be chosen later but satisfy the property that {
		$\cM:=\cF\cF^* + \tilde{\cS} \succ 0$ and 
		$\cN:=\cG\cG^*+\tilde{\cT} \succ 0$.} For $k=0,1,...$ perform the following steps:
		
		\textbf{Step 1.} Compute 
		\begin{eqnarray*}
			&&{\mu}^{k+1}\approx {\rm argmin}_{\mu\in (\bcX\times \cX)\times \cY } \big\{\cL_\sigma (\mu,\nu^k;\omega^k)+\frac{\sig}{2}\norm{\mu-\mu^k}_{\tilde{\cS}}^2\big\}, \\
		\end{eqnarray*}
		\blue{
		such that there exists a vector $\delta^k$ satisfying $\|\cM^{-\frac{1}{2}}{\delta}^k\| \le \tilde{\eps}_k$ and
		\[
			{\delta}^k \in \partial_{\mu}\cL_\sigma (\mu^{k+1},\nu^k;\omega^k) + \sigma\tilde{\cS}({\mu}^{k+1}-\mu^k).
		\]
		}
		\textbf{Step 2.} Compute 
		\begin{eqnarray*}
			&&\nu^{k+1}\approx \mbox{argmin}_{\nu\in (\cX\times\bcX)\times\bcY}
			\big\{\cL_\sigma (\mu^{k+1},\nu;\omega^k)
			+\frac{\sig}{2}\norm{\nu-\nu^k}_{\tilde{\cT}}^2\big\}, \\
		\end{eqnarray*}
		\blue{
		such that there exists a vector $\delta^k$ satisfying $\|\cN^{-\frac{1}{2}}{\gamma}^k\| \le \tilde{\eps}_k$ and 
		\[
			\gamma^k \in \partial_{\nu}\cL_\sigma (\mu^{k+1},\nu^{k+1};\omega^k) +\sigma \tilde{\cT}(\nu^{k+1}-\nu^k).
		\]
		}
		\textbf{Step 3.} Compute $\omega^{k+1}:=\omega^k+\tau\sigma (\cF^* \mu^{k+1} + \cG^* \nu^{k+1} - h)$.
}}

\bigskip
The convergence of the above imsPADMM algorithm is proven in
\cite{CST}.
\section{Inexact symmetric Gauss-Seidel proximal ADMM and ALM for solving (D)}

We now apply the abstract inexact proximal ADMM (imsPADMM
presented in the last section)
to the DBA problem. Explicitly, given a penalty parameter $\sig > 0$, 
the augmented Lagrangian function associated with (D) is given as follows:
for 
$(\bz,z,y)\in \bcX \times \cX \times \cY $,
$(v,\bar{v},\by)\in \cX\times\bcX\times\bcY$,
$(x,\bx)\in\cX\times\bcX$,
\begin{eqnarray}
&& \hspace{-7mm}
\cL_\sig(\bz,z,y;v,\bar{v},\by;x,\bx) \;=\; \theta^*(-v) +\delta_\cK^*(-z) -\inprod{b}{y} + \btheta^*(-\bar{v})
+\delta_{\bcK}^*(-\bz) -\inprod{\bb}{\by} \nn \\[5pt]
&& \hspace{3cm} 
-\frac{1}{2\sig}\norm{x}^2 -\frac{1}{2\sig}\norm{\bx}^2 +\frac{\sigma}{2}\norm{A^*y+B^*\by+z + v -c + \sig^{-1}x}^2 
\nn \\[5pt]
&&  \hspace{3cm}
+\frac{\sigma}{2}\norm{\bB^*\by+\bz+ \bar{v} -\bc+\sig^{-1}\bx}^2 .
\end{eqnarray}


Next we present the explicit
algorithmic framework of our
inexact sGS proximal ADMM  
for solving the dual problem (D)  of the block-angular 
problem (P). To apply the sGS decomposition theorem within the inexact proximal ADMM {framework}, 
we consider $(\bz,z,y)$ as the first group of block variables and $(v,\bv,\by)$ as the second group
of block variables. 
Then we make use of the sGS decomposition theorem in Step 1 of imsPADMM
by considering the quadratic term $\cF\cF^*$ corresponding to the block $\mu$
in $\cL_\sig(\mu,\nu^k;\omega^k)$.
Similarly, we make {use} of the sGS decomposition in Step 2 of imsPADMM
by considering $\cG\cG^*$ corresponding to the block $\nu$
in $\cL_\sig(\mu^{k+1},\nu;\omega^k)$.

The explicit template of the inexact sGS proximal ADMM for solving (D) is given as follows.

\bigskip
\noindent{{\bf Algorithm DBA-sGS-ADMM}.} 
Let $\{\eps_k\}$ be a given nonnegative summable sequence, {i.e. $\sum_{k=0}^\infty\epsilon_k<\infty$}.
Given the initial points
$(\bz^0,z^0,y^0)\in \bcX\times \cX\times\cY$, $(v,\bar{v}^0,\by^0) \in\cX\times\bcX\times\bcY$, $(x^0,\bx^0)\in\cX\times\bcX$, perform the following steps in each iteration. 
\begin{description}
	\item[Step 1.]  Let $J:\cY\to \cY$ be a given self-adjoint positive semidefinite linear operator, 
	and $\cS$ be a positive semidefinite linear operator that is to be chosen later. 
	Compute
	\begin{eqnarray*}
		(\bz^{k+1},z^{k+1},y^{k+1}) \approx \argmin\limits_{\bz,z,y} \left\{
		\begin{array}{l}
			\cL_\sig (\bz,z,y;v^k,\bar{v}^k,\by^k;x^k,\bx^k)
			\;+
			\\[5pt]
			\frac{\sig}{2} \norm{y-y^k}^2_J
			+ \frac{\sig}{2}\norm{(z;y)-(z^k;y^k)}_\cS^2
		\end{array}
		\Bigm\vert \begin{array}{l}
			y\in\cY, 
			\\[3pt] 
			z\in\cX,\bz\in\bcX
		\end{array}
		\right\}.
	\end{eqnarray*}
	Next we explain how the above abstract subproblem can be solved in practice.
	Let $R^{1,k} = A^*y^k+B^*\by^k+z^k+v^k  -c^k$,
	$R^{2,k} = \bB^*\by^k+\bz^k +\bar{v}^k -\bc^k$, where $c^k =c - \sig^{-1}x^k$
	and $\bc^k = \bc - \sig^{-1}\bx^k.$
	Observe that $\bz$ and $(z,y)$ are separable. Hence we can compute $\bz^{k+1}$ and $(z^{k+1},y^{k+1}) $
	separately by solving
	\begin{eqnarray}
	\bz^{k+1} &=& \argmin\limits_{\bz}  \Big\{
	\frac{\sigma}{2}\norm{\bz -\bz^k + R^{2,k}}^2  + \delta_{\bcK}^*(-\bz) \bigm\vert \bz \in \bcX
	\Big\}
	\\[5pt]
	(z^{k+1},y^{k+1}) & \approx & 
	\argmin\limits_{z,y} \left\{
	\begin{array}{l|}
	-\inprod{b}{y} + \delta_\cK^*(-z) 
	\\[5pt]
	+\frac{\sigma}{2}\norm{z+A^*y  + B^*\by^k+v^k-c^k}^2  
	\\[5pt]
	\frac{\sig}{2} \norm{y-y^k}^2_J
	+ \frac{\sig}{2}\norm{(z;y)-(z^k;y^k)}_\cS^2
	\end{array}
	\;  \begin{array}{l} y\in\cY, \\[3pt] z \in \cX \end{array}
	\right\}.
	\qquad \label{eq-yz}
	\end{eqnarray}
	Notice that $\bz^{k+1}$ can be computed by solving the following $N$ subproblems in {\em parallel}, i.e.,
	\begin{eqnarray}
	\bz_i^{k+1} &=& 
	\argmin\limits_{\bz_i}  \Big\{
	\frac{\sigma}{2}\norm{\bz_i -\bz^k_i + R^{2,k}_i}^2  + \delta_{\cK_i}^*(-\bz_i) \bigm\vert \bz_i \in \cX_i
	\Big\}
	\\[5pt]
	&=& 
	-{\rm Prox}_{\sig^{-1} \delta_{\cK_i}^*}\big( R^{2,k}_i-\bz^k_i\big), \quad i=1,\ldots,N. 
	\nonumber
	\end{eqnarray}
	{In order to compute $(z^{k+1},y^{k+1})$ efficiently in a sequential manner, we apply the sGS decomposition
	theorem to the problem \eqref{eq-yz}. We regard $z$ and $y$ in \eqref{eq-yz} as $x_1$ and $x_2$ in the sGS decomposition
	theorem, respectively.} We then choose
	$\cS = {\bf sGS}(\cQ_{z,y})$ to be the sGS linear operator  
	\eqref{eq-sGS} associated with the following positive semidefinite linear operator
	to decompose the computation of $y^{k+1}$ and 
	$z^{k+1}$: 
	\begin{eqnarray*}
		\cQ_{z,y} : = \left[\begin{array}{cc} I_n & A^* \\[3pt]  A & AA^*+J
		\end{array}\right]. 
	\end{eqnarray*}
	{In the context of imsPADMM, we have chosen $\tilde{\cS} = {\rm diag}(0, {\rm diag}(0,J)+{\bf sGS}(\cQ_{z,y})).$}
	Now we can compute
	$(z^{k+1},y^{k+1})$ in a symmetric Gauss-Seidel fashion as follows:
	\begin{eqnarray*}
		\begin{array}{lll}
			y^{k+1}_{\rm tmp}& = & \argmin\limits_{y} 
			\big\{-\inprod{b}{y}
			+\frac{\sigma}{2}\norm{A^*(y-y^k)+R^{1,k}}^2  + \frac{\sig}{2}\norm{y-y^k}_J^2 
			-\inprod{\delta^{k}_{\rm tmp}}{y}
			\big\}
			\\[8pt]
			& =& y^k + (AA^* + J )^{-1} \big(\sig^{-1}(b +\delta^{k}_{\rm tmp}) - AR^{1,k}\big)
			\\[8pt]
			z^{k+1} &=& \argmin\limits_{z} \big\{
			\frac{\sigma}{2}\norm{z-z^k+A^*(y^{k+1}_{\rm tmp} -y^k)+R^{1,k}}^2 + \delta_\cK^*(-z) \mid z \in \cX
			\big\}
			\\[8pt]
			&=&  - {\rm Prox}_{\sig^{-1}\delta_\cK^*} \big( A^*(y^{k+1}_{\rm tmp} -y^k)+R^{1,k}-z^k\big)
			\\[8pt]
			y^{k+1} &= & \argmin\limits_{y}  \big\{-\inprod{b}{y}
			+\frac{\sigma}{2}\norm{ A^*(y-y^k)+R^{1,k}+ z^{k+1}-z^k }^2 
			+ \frac{\sig}{2}\norm{y-y^k}_J^2 -\inprod{\delta^{k}}{y}
			\big\}
			\\[5pt]
			&=& y^{k} + (AA^*+J)^{-1} (\sig^{-1}(b+\delta^{k}) - A(R^{1,k} + z^{k+1}-z^k)) 
		\end{array}
	\end{eqnarray*}
	where the error vectors 
	${\delta}^{k}_{\rm tmp}$ and 
	${\delta}^{k}$ satisfy the condition that
	$\max\{\norm{{\delta}^{k}_{\rm tmp}},\norm{\delta^{k}}\} \leq \eps_k$.

	\item[Step 2.] Let $\bJ$ be a given block diagonal
	linear operator with $\bJ =\diag{\bJ_1,\ldots,\bJ_N}$ such that each $\bJ_i :\cY_i\to \cY_i$ is 
	self-adjoint positive semidefinite, 
	and \blue{$\cT:\cX\times\bar\cX\times\bar\cY\to\cX\times\bar\cX\times\bar\cY$} be a  positive semidefinite linear operator that is to be chosen later. 
	Compute 
	\begin{eqnarray*}
		(v^{k+1},\bar{v}^{k+1},\by^{k+1}) \approx 
		\argmin\limits_{v,\bv,\by} 
		\left\{ 
		\begin{array}{l}
			\cL_\sig (\bz^{k+1},z^{k+1},y^{k+1};v,\bar{v},\by;x^k,\bx^k)
			\; +
			\\[5pt]
			\frac{\sig}{2}\norm{\by-\by^k}_{\bJ}^2
			+ \frac{\sig}{2}\norm{( v-v^k;\bar{v}-\bar{v}^k; \by-\by^k)}_\cT^2
		\end{array}
		\Bigm\vert  \begin{array}{l}
			\by\in\RR^{\bar{m}},\\
			v\in\cX,\\
			\bar{v}\in\bcX
		\end{array}
		\right\}.
	\end{eqnarray*}
	Next we explain how the above abstract subproblem can be solved efficiently in practice.
	Let 
	\begin{eqnarray}
	M = BB^* + \bB\bB^* + \bJ ,
	\label{eq-M}
	\end{eqnarray}
	and 
	$R^{3,k} = A^*y^{k+1}+B^*\by^k+z^{k+1}+v^k  -c^k = R^{1,k} + A^*(y^{k+1}-y^k)+ (z^{k+1}-z^k)$,
	$R^{4,k} = \bB^*\by^k+\bz^{k+1} +\bar{v}^k -\bc^k = R^{2,k} + (\bz^{k+1}-\bz^k)$. 
	{In order to compute $(v^{k+1},\bar{v}^{k+1},\by^{k+1})$ efficiently in a sequential manner, 
	we regard $(v,\bar v)$ and $\bar y$ in the subproblem as $x_1$ and $x_2$ in the sGS decomposition theorem, respectively. We then choose  $\cT = {\bf sGS}(\cQ_{v\bv,\by})$, which is the sGS linear operator associated with the following positive semidefinite linear operator:}
	\begin{eqnarray*}
		\cQ_{v\bv,\by}  = \left[ 
		\begin{tabular}{cc|c}
			$I_n$ & 0  & $B^*$ 
			\\{}
			0 &$I_{\bn}$ & \small{$\bB^*$}
			\\ \hline
			$B$ &\blue{\small{$\bB$}}  & $M$
		\end{tabular}
		\right].
	\end{eqnarray*}
	 {In the context of imsPADMM, we have chosen $\tilde{\cT} = {\rm diag}(0,0,\bar{J}) + {\bf sGS}(\cQ_{v\bv,\by}).$}
	Then we can decompose the computation of $(v^{k+1},\bar{v}^{k+1},\by^{k+1}) $
	in a symmetric Gauss-Seidel fashion
	as follows:
	\begin{eqnarray}
	\label{eq-step2}
	\begin{array}{l}
	\by^{k+1}_{\rm tmp}\;=\; \argmin\limits_{\by} 
	\left\{
	\begin{array}{l}
	-\inprod{\bb}{\by}
	+\frac{\sigma}{2}\norm{B^*(\by-\by^k) +R^{3,k}}^2 
	+\frac{\sigma}{2}\norm{\bB^*(\by-\by^k)+ R^{4,k}}^2 
	\\[5pt]
	+\frac{\sig}{2}\norm{\by-\by^k}_{\bJ}^2 -\inprod{{\bar{\delta}^{k}_{\rm tmp}}}{\by}
	\end{array}
	\right\}
	\\[20pt]
	\qquad\;\; = \;
	\by^k + M^{-1} (\sig^{-1}(\bb+\bar{\delta}^{k}_{\rm tmp}) - B R^{3,k} -\bB R^{4,k})
	\\[15pt]
	\left\{ 
	\begin{array}{rcl}
	v^{k+1} &=& \argmin\limits_{v}  \big\{ \theta^*(-v) +
	\frac{\sigma}{2}\norm{v-v^k + B^*(\by^{k+1}_{\rm tmp}-\by^k) +R^{3,k}}^2
	\bigm\vert  v\in\cX
	\big\}
	\\[8pt]
	\bar{v}^{k+1} &=& \argmin\limits_{\bv}  \big\{ \btheta^*(-\bar{v})  + 
	\frac{\sigma}{2}\norm{\bar{v}-\bar{v}^k + \bB^*(\by^{k+1}_{\rm tmp}-\by^k) +R^{4,k}}^2 
	\bigm\vert \bar{v}\in\bcX
	\big\}
	\end{array} \right.
	\\[25pt]
	\by^{k+1} \;  =  \; \argmin\limits_{\by}  \left\{ 
	\begin{array}{l} 
	-\inprod{\bb}{\by}
	+\frac{\sigma}{2}\norm{B^*(\by-\by^k)+v^{k+1}-v^k+R^{3,k}}^2 
	\\[5pt]
	+\frac{\sigma}{2}\norm{\bB^*(\by-\by^k)+\bar{v}^{k+1}-\bar{v}^k+ R^{4,k} }^2
	+\frac{\sig}{2}\norm{\by-\by^k}_{\bJ}^2 - \inprod{{\bar{\delta}^{k}}}{\by}
	\end{array}
	\right\}
	\\[20pt]
	\qquad\;\; = \; \by^{k} + M^{-1}\big(\sig^{-1}(\bb+\bar{\delta}^{k}) - B (R^{3,k} +v^{k+1}-v^k) 
	-\bB (R^{4,k} + \bar{v}^{k+1}-\bar{v}^k)\big) 
	\end{array}
	\end{eqnarray}
	where the error vectors 
	${\bar{\delta}}^{k}_{\rm tmp}$ and 
	${\bar{\delta}}^{k}$ satisfy the condition that
	$\max\{\norm{{\bar{\delta}}^{k}_{\rm tmp}},\norm{\bar{\delta}^{k}}\} \leq \eps_k$.
	
	We should note that the computation of $\bar{v}^{k+1}$ in \eqref{eq-step2} can be done
	by solving $N$ smaller subproblems  in 
	{\em parallel}, where 
	\begin{eqnarray}
	\bar{v}^{k+1}_i &=&
	\mbox{argmin} \Big\{ \btheta_i^*(-\bar{v}_i) + 
	\frac{\sigma}{2}\norm{\bar{v}_i-\bar{v}_i^k + \bB_i^*\big((\by^{k+1}_{\rm tmp})_i -\by_i^k\big) +{R^{4,k}_i}}^2 
	\bigm\vert \bar{v}_i \in \cX_i
	\Big\} 
	\nonumber
	\\[5pt]
	&= & - {\rm Prox}_{\sig^{-1}\btheta^*_i} \Big( \bB_i^*\big((\by^{k+1}_{\rm tmp})_i -\by_i^k\big) +R^{4,k}_i -\bar{v}_i^k\Big), \quad i=1,\ldots,N.
	\end{eqnarray}

	\item[Step 3.] Compute
	\begin{eqnarray*}
		x^{k+1} &=& x^k + \tau \sig (A^* y^{k+1} + B^*\by^{k+1}+z^{k+1}+v^{k+1}-c),
		\\[5pt]
		\bx^{k+1} &=& \bx^k + \tau \sig (\bB^*\by^{k+1}+\bz^{k+1}+\bar{v}^{k+1}-\bc),
	\end{eqnarray*}
	where $\tau \in (0,(1+\sqrt{5})/2)$ is the step-length that is usually set
	to be $1.618$. 
\end{description}

\blue{
Let $\cU_\mu$ and $\cD_\mu$ be the upper block triangular part and block diagonal part of $\cQ_{z,y}$ respectively, $\cU_\nu$ and $\cD_\nu$ be the upper block triangular part and block diagonal part of $\cQ_{v\bar v,\bar y}$ respectively. Define 
\begin{equation}
\begin{aligned}
&\Delta^k_\mu=\left[\begin{array}{c}0\\\widetilde\Delta^k_\mu\end{array}\right],\;\widetilde\Delta^k_\mu=\left[\begin{array}{c}0\\\delta^k\end{array}\right]+\cU_\mu\cD_\mu^{-1}\left(\left[\begin{array}{c}0\\\delta^k\end{array}\right]-\left[\begin{array}{c}0\\\delta_{\rm tmp}^k\end{array}\right]\right),\\
&\Delta_\nu^k=\left[\begin{array}{c}0\\\delta^k\end{array}\right]+\cU_\nu\cD_\nu^{-1}\left(\left[\begin{array}{c}0\\\bar\delta^k\end{array}\right]-\left[\begin{array}{c}0\\\bar\delta_{\rm tmp}^k\end{array}\right]\right).
\end{aligned}
\label{Delta-ADMM-def}
\end{equation}
Based on the above notations and the
practical implementations of Step 1 and Step 2 in DBA-sGS-ADMM, we have the following results, which is a direct consequence of Theorem \ref{sGS-decom}.

\begin{lemma}
Using the notations in Section 4.2, let $\tilde{\cS} = {\rm diag}(0, {\rm diag}(0,J)+{\bf sGS}(\cQ_{z,y}))$, $\tilde{\cT} = {\rm diag}(0,0,\bar{J}) + {\bf sGS}(\cQ_{v\bv,\by})$, suppose $\cD_\mu\succ0$, $\cD_\nu\succ0$. The iterates $\{(\mu^k,\nu^k)\}$ generated by DBA-sGS-ADMM is the unique solution to the following perturbed proximal minimization problem:
\[
\begin{aligned}
\mu^{k+1}&=\argmin_{\mu\in\bar\cX\times\cX\times\cY}\left\{\cL_{\sigma}(\mu,\nu^k;\omega^k)+\frac{\sigma}{2}\|\mu-\mu^k\|^2_{\tilde\cS}-\inprod{\Delta_\mu^k}{\mu}\right\}\\
\nu^{k+1}&=\argmin_{\nu\in\cX\times\bar\cX\times\bar\cY}\left\{\cL_{\sigma}(\mu^{k+1},\nu;\omega^k)+\frac{\sigma}{2}\|\nu-\nu^k\|^2_{\tilde\cT}-\inprod{\Delta_\nu^k}{\nu}\right\},
\end{aligned}
\]
which is equivalent to the following inclusion problem:
\[
\begin{aligned}
\Delta_\mu^{k}&\in\partial_\mu\cL_{\sigma}(\mu,\nu^k;\omega^k)+\sigma\tilde\cS(\mu-\mu^k)\\
\Delta_\nu^{k}&\in\partial_\nu\cL_{\sigma}(\mu^{k+1},\nu;\omega^k)+\sigma\tilde\cT(\nu-\nu^k).
\end{aligned}
\]
\label{sGS-ADMM-equi-lemma}
\end{lemma}
}
\blue{From the inexact conditions in DBA-sGS-ADMM and \eqref{Delta-ADMM-def}, we can obtain that $\|\Delta^k_\mu\|$ and $\|\Delta_\nu^k\|$ are summable as long as $\{\epsilon_k\}$ is summable.
}

\subsection{An inexact symmetric Gauss-Seidel proximal ALM for solving (D)}

When the functions $\theta(\cdot)$ and $\btheta_i(\cdot)$ ($i=1,\ldots,N$) are absent in 
\eqref{eq-DBA}, we can design a symmetric Gauss-Seidel proximal augmented Lagrangian method (sGS-ALM) to
solve the problem. In this case, the variables $v$ and $\bar{v}$ are absent in (D), and the 
associated augmented Lagrangian function is given as follows: 
for 
$\blue{(z,\bz,y,\by)}\in \cX\times\bcX\times \cY\times \bcY$, 
$(x,\bx)\in\cX\times\bcX$,
\begin{eqnarray}
\cL_\sig(z,\bz,y,\by;x,\bx) &=& -\inprod{b}{y} -\inprod{\bb}{\by} + \delta_\cK^*(-z) + \delta_\bcK^*(-\bz) -\frac{1}{2\sig}\norm{x}^2 -\frac{1}{2\sig}\norm{\bx}^2 \nn \\[5pt]
&&
+\frac{\sigma}{2}\norm{A^*y+B^*\by+z -c + \sig^{-1}x}^2 
\nn \\[5pt]
&&
+\frac{\sigma}{2}\norm{\bB^*\by+\bz -\bc+\sig^{-1}\bx}^2 .
\end{eqnarray}

\noindent{{\bf Algorithm DBA-sGS-ALM}.} Let $\{\eps_k\}$ be a given nonnegative summable sequence, {i.e. $\sum_{k=0}^\infty\epsilon_k<\infty$}. Given the initial points
$(y^0,\by^0,z^0,\bz^0) \in \cY\times \bcY\times\cX\times\bcX$, $(x^0,\bx^0)\in\cX\times\bcX$, perform the following steps in each iteration. 
\begin{description}
	\item[Step 1.]  Let $J:\cY\to\cY$ and $\bar{J}:\bcY \to \bcY$ 
	be given self-adjoint positive semidefinite linear operators, 
	and $\cS$ be a positive semidefinite linear operator that is to be chosen later. 
	Compute
	\begin{eqnarray*}
		(\bz^{k+1},z^{k+1},y^{k+1},\by^{k+1}) \approx  
		\mbox{argmin} 
		\left\{
		\begin{array}{l|}
			\cL_\sig (z,\bz,y,\by;x^k,\bx^k) \;
			\\[5pt]
			+ \frac{\sig}{2} \norm{y-y^k}^2_J  + \frac{\sig}{2} \norm{\by-\by^k}^2_{\bar{J}}
			\\[5pt]
			+ \frac{\sig}{2}\norm{(\bz;z;y;\by)-(\bz^k;z^k; y^k;\by^k)}_\cS^2
		\end{array}
		\; \; 
		\begin{array}{l}
			y\in\cY, \\ \by\in \bcY, 
			\\
			z\in\cX,\\ \bz\in\bcX
		\end{array} 
		\right\}.
	\end{eqnarray*}
	The inexact conditions will be clarified later. 
	Note that the quadratic term of the above objective function 
	involves the following positive semidefinte linear operator
	\begin{eqnarray*}
		\cQ_{\bz z, y, \by} = \left[ 
		\begin{tabular}{cc|c|c} 
			$I_{\bn}$   & 0 & 0  & $\bB^*$
			\\ 
			0 & $I_{n}$ &  $A^*$  & $B^*$ 
			\\ \hline
			0 & $A$  & $AA^* + J$   &$AB^*$
			\\ \hline
			\blue{\small{$\bB$}} & $B$  &$BA^*$ &\small{$BB^* + \blue{\bB\bB^*} + \bar{J}$}
		\end{tabular}
		\right].
	\end{eqnarray*}
	{In order to compute $(\bar z^{k+1},z^{k+1},y^{k+1},\bar y^{k+1})$ efficiently in a sequential manner,  we regard $(\bar z,z)$, $y$ and $\bar y$ in the subproblem as $x_1$, $x_2$ and $x_3$ in the sGS decomposition theorem, respectively.} By choosing $\cS = {\bf sGS}(\cQ_{\bz z, y,\by})$, the sGS linear
	operator associated with $\cQ_{\bz z, y, \by}$,
	then we can decompose the computation of $(\bz^{k+1},z^{k+1},y^{k+1},\by^{k+1})$
	in a symmetric Gauss-Seidel fashion as follows.
	Let  $c^k = c-\sig^{-1} x^k$ and 
	$\bc^k = \bc - \sig^{-1}\bx^k.$ Compute
	\begin{eqnarray*}
		&&\by^{k+1}_{\rm tmp} \; = \; \mbox{argmin}\left\{
		\begin{array}{l}
			-\inprod{\bb}{\by}
			+\frac{\sigma}{2}\norm{B^*\by  + A^* y^{k} + z^k - c^k}^2 
			+\frac{\sigma}{2}\norm{\bB^*\by + \bz^k - \bc^k}^2 
			\\[3pt]
			+\frac{\sig}{2}\norm{\by-\by^k}_{\bJ}^2 - \inprod{\bar{\delta}^{k}_{\rm tmp}}{\by}
		\end{array}
		\right\}
		\\[8pt]
		&& y^{k+1}_{\rm tmp} \; = \; \mbox{argmin}\big\{-\inprod{b}{y}
		+\frac{\sigma}{2}\norm{A^*y + B^*\by^{k+1}_{\rm tmp} +z^k-c^k}^2  + \frac{\sig}{2}\norm{y-y^k}_J^2
		-\inprod{\delta^{k}_{\rm tmp}}{y}
		\big\}
		\\[8pt]
		& &\left\{  \begin{array}{rcl}
			z^{k+1} &=&
			\mbox{argmin}\big\{
			\frac{\sigma}{2}\norm{z+A^*y^{k+1}_{\rm tmp} + B^*\by_{\rm tmp}^{k+1}-c^{k}}^2 + \delta_\cK^*(-z) 
			\mid z \in \cX
			\big\}
			\\[8pt]
			&=&
			-{\rm Prox}_{\sig^{-1}\delta_\cK^*}\big( A^*y^{k+1}_{\rm tmp} + B^*\by_{\rm tmp}^{k+1}-c^{k}\big)
			\\[8pt]
			\bz^{k+1}_i &=&
			\mbox{argmin} \big\{
			\frac{\sigma}{2}\norm{\bz_i+ \bB_i^* (\by_{\rm tmp}^{k+1})_i -\bc^k_i}^2 +\delta_{\cK_i}^*(-\bz_i)  
			\mid \bz_i \in \cX_i
			\big\}
			\\[5pt]
			&=& -{\rm Prox}_{\sig^{-1}\delta_{\cK_i}^*} \big( \bB_i^* (\by_{\rm tmp}^{k+1})_i -\bc^k_i
			\big)
			\quad i=1,\ldots,N
		\end{array}\right.
		\\[8pt]
		&& y^{k+1} \; = \;\mbox{argmin}\big\{-\inprod{b}{y}
		+\frac{\sigma}{2}\norm{ A^*y + B^*\by^{k+1}_{\rm tmp}+ z^{k+1}-c^k }^2 
		+ \frac{\sig}{2}\norm{y-y^k}_J^2 - \inprod{\delta^{k}}{y}
		\big\}
		\\[8pt]
		&& \by^{k+1}\; = \; \mbox{argmin}\left\{
		\begin{array}{l}
			-\inprod{\bb}{\by}
			+\frac{\sigma}{2}\norm{B^*\by+A^*y^{k+1}+ z^{k+1}-c^k}^2			
			\\[5pt]
			+\frac{\sigma}{2}\norm{\bB^*\by+ \bz^{k+1}-\bc^k}^2 
			+\frac{\sig}{2}\norm{\by-\by^k}_{\bJ}^2 -\inprod{\bar{\delta}^{k}}{\by}
		\end{array}
		\right\}
	\end{eqnarray*}
	where the error vectors $\delta_{\rm tmp}^k$, $\bar{\delta}_{\rm tmp}^k$,
	$\delta^k$, $\bar{\delta}^k$ satisfy
	$\max\{\norm{\delta_{\rm tmp}^k},\norm{\bar{\delta}_{\rm tmp}^k},
	\norm{\delta^k},\norm{\bar{\delta}^k}\}\leq \eps_k$.
	
	\item[Step 2.] Compute
	\begin{eqnarray*}
		x^{k+1} &=& x^k + \tau \sig (A^* y^{k+1} + B^*\by^{k+1}+z^{k+1}-c),
		\\[5pt]
		\bx^{k+1} &=& \bx^k + \tau \sig (\bB^*\by^{k+1}+\bz^{k+1}-\bc),
	\end{eqnarray*}
	where $\tau \in (0,2)$ is the step-length which is usually set
	to be $1.9$. 
\end{description}

\blue{
Let $\cU$ and $\cD$ be the upper block triangular part and the block diagonal part of $\cQ_{\bar z z,y,\bar y}$, respectively, and
\begin{equation}
\Delta^k:=\left[\begin{array}{c} 0\\\delta^k\\\bar\delta^k
		\end{array}\right]+\cU\cD^{-1}\left(\left[\begin{array}{c} 0\\\delta^k\\\bar\delta^k
		\end{array}\right]-\left[\begin{array}{c} 0\\\delta_{\rm tmp}^k\\\bar\delta_{\rm tmp}^k
		\end{array}\right]\right).
		\label{Delta-def}
\end{equation}
The summable condition for $\epsilon_k$ and the inexact condition: $\max\{\norm{\delta_{\rm tmp}^k},\norm{\bar{\delta}_{\rm tmp}^k},
	\norm{\delta^k},\norm{\bar{\delta}^k}\}\leq \eps_k$ in Step 1 can guarantee that the error $\|\Delta^k\|$ is also summable, since
\[
\|\Delta^k\|\leq\|\delta^k\|+\|\bar \delta^k\|+\|\cU\|\|\cD^{-1}\|\left(\|\delta^k\|+\|\bar\delta^k\|+\|\delta^k_{\rm tmp}\|+\|\bar\delta^k_{\rm tmp}\|\right),
\]
{where the norms $\|\cU\|$ and $\|\cD^{-1}\|$ here are induced operator norms defined in Section 1.}
Based on the above notations and the practical implementation of Step 1, we have the following results, which is a direct consequence of Theorem \ref{sGS-decom}.

\begin{lemma}
Suppose $\cD\succ0$, let $w^{k+1}:=(z^{k+1},\bar z^{k+1},y^{k+1},\bar y^{k+1})$ and $\tilde\cS=\diag{0,0,J,\bar J}+ {\bf sGS}(\cQ_{\bz z, y,\by})$. The iterates $w^{k+1}$ generated by Step 1 in DBA-sGS-ALM is the unique solution to the following perturbed proximal minimization problem:
\[
w^{k+1}=\argmin_{w\in\cX\times\bar\cX\times\cY\times\bar\cY}\left\{\cL_\sigma(w;x^k,\bar x^k)+\frac{\sigma}{2}\|w-w^k\|_{\tilde\cS}^2-\inprod{\Delta^k}{w}\right\},
\]
which is equivalent to the following inclusion problem
\[
\Delta^k\in\partial_w\cL_\sigma(w;x^k,\bar x^k)+\sigma\tilde\cS(w-w^k).
\]
\label{sGS-ALM-equi-lemma}
\end{lemma}
}

\subsection{Convergence of the inexact DBA-sGS-ADMM and DBA-sGS-ALM algorithms}

{The convergence theorem of the inexact DBA-sGS-ADMM and DBA-sGS-ALM algorithms could be established directly using known results from \cite{CST}, \cite{CLST} and \cite{LW-imiPADMM}.} Here we would state the theorems again for the convenience of the readers. 

In order to understand the theorem, first we need some definitions.
\begin{definition}
	Let $\mathbb{F}:\cX \rightrightarrows\cY$ be a multivalued mapping and denote its inverse by $\mathbb{F}^{-1}$. The graph of the multivalued function $\mathbb{F}$ is defined by 
	${\rm gph}\mathbb{F}:=\{(x,y)\in\cX\times \cY\mid y\in \mathbb{F}(x)\}$.
\end{definition}

Denote $u:=(y,\by,v,\bv,z,\bz,x,\bx)\in\mathcal{U}:=\cY\times\bcY\times\cX\times\bcX\times\cX\times\bcX\times\cX\times\bcX$. The KKT mapping $\mathcal{R}:\mathcal{U}\rightarrow\mathcal{U}$ of \eqref{eq-DBA} is defined by 
\begin{eqnarray}
\mathcal{R}(u):=\begin{pmatrix}
Ax - b \\
Bx + \bB\bx - \bb \\
A^*y + B^* \by + z + v - c \\
\bB^*\by + \bz + \bar{v} - \bc \\
x - \Pi_{\cK}(x-z)  \\
\bx - \Pi_{\bcK}(\bx-\bz) \\
x - {\rm Prox}_{\theta}(x-v) \\
\bx - {\rm Prox}_{\btheta}(\bx-\bar{v})
\end{pmatrix}.
\end{eqnarray}

Denote the set of KKT points by $\bar{\Omega}$. The  KKT mapping $\mathcal{R}$ is said to be metrically subregular at $(\bar{u},0)\in {\rm gph}\mathcal{R}$ with modulus $\eta>0$ if
there exists a scalar $\rho>0$ such that 
$${\rm dist}(u,\bar{\Omega})\le \eta\|\mathcal{R}(u)\| \; \; \forall u \in \{u\in \mathcal{U}: \|u-\bar{u}\|\le \rho\}.$$

Now we are ready to present the convergence theorem of DBA-sGS-ADMM.
\begin{theorem}
	\blue{Let $\Delta_\mu^k$ and $\Delta_\nu^k$ be defined as \eqref{Delta-ADMM-def}, $\cM:=\cF\cF^* + \tilde{\cS} \succ 0$ and 
		$\cN:=\cG\cG^*+\tilde{\cT} \succ 0$,}
	 $\{u^k:=(y^k,\by^k,v^k,\bv^k,z^k,\bz^k,x^k,\bx^k)\}$ be the sequence generated by DBA-sGS-ADMM, where $\tilde{\cS} = {\rm diag}(0, {\rm diag}(0,J)+{\bf sGS}(\cQ_{z,y}))$, $\tilde{\cT} = {\rm diag}(0,0,\bar{J}) + {\bf sGS}(\cQ_{v\bv,\by})$. Assume the KKT point set $\bar \Omega$ is nonempty. Then, we have the following results.
	
	\begin{enumerate}
		\item[(a)] The sequence $\{(y^k,\by^k,v^k,\bv^k,z^k,\bz^k)\}$ converges to an optimal solution of the dual problem (D) of \eqref{eq-DBA}
		and the sequence $\{(x^k,\bx^k)\}$ converges to an optimal solution of the  primal problem (P) of \eqref{eq-DBA}.
		
		\blue{\item[(b)] Suppose that the sequence $\{u^k\}$ converges to a KKT point $\hat{u}:=(\hat{y},\hat{\by},\hat{v},\hat{\bv},\hat{z},\hat{\bz},\hat{x},\hat{\bx})$ and the KKT mapping $\mathcal{R}$ is metrically subregular at $(\hat{u},0)\in {\rm gph}\mathcal{R}$. If a relative error criterion is employed instead:
		\[
		\begin{aligned}
			\|\Delta_\mu^k\|_{\cM^{-1}}&\leq\eta_k\left(\|\mu^{k}-\mu^{k+1}\|^2_\cM+\|\mu^{k-1}-\mu^k\|_\cM^2\right)^{\frac{1}{2}},\\
			\|\Delta_\nu^k\|_{\cN^{-1}}&\leq\eta_k\left(\|\nu^{k}-\nu^{k+1}\|^2_\cN+\|\nu^{k-1}-\nu^k\|_\cN^2\right)^{\frac{1}{2}}.
		\end{aligned}
		\]
		and $\eta_k\in(0,1)$, $\sum_{k\geq0}\eta_k^2<\infty$, then the sequence $\{u^k\}$ is linearly convergent to $\hat{u}$.}
	\end{enumerate}
	
\end{theorem}

\begin{Proof} 
	The result of (a) follows directly by applying the convergence result in {\cite[Theorem 5.1]{CST} to (D). The result of (b) follows directly from \cite[Theorem 1]{LW-imiPADMM}}.
\end{Proof}

\begin{remark} By Theorem 1 and Remark 1 in \cite{LSTc}, we know that when (P) is 
	a convex programming problem where for each $i=1,\ldots,N$,  $\theta$ and $\btheta_i$ are piecewise linear-quadratic or strongly convex, and
	$\cK,\cK_i$ are polyhedral, then 
	$\cal R$ is metrically 
	subregular at $(\bar{u},0)\in {\rm gph}\mathcal{R}$  for any KKT point $\bar{u}.$
\end{remark}

\blue{Implied by Lemma \ref{sGS-ALM-equi-lemma}, DBA-sGS-ALM is equivalent to an inexact majorized indefinite-proximal ALM. Therefore, the convergence of the DBA-sGS-ALM can be established readily by using known
results in \cite[Theorem 3.1 and Theorem 3.2]{CLST}.} We state the convergence theorem here 
for the convenience of the readers.

\begin{theorem}
	\blue{Let $\Delta^k$ be defined as \eqref{Delta-def}, $\{\eta^k\}$ be a given sequence of nonnegative numbers that converges to 0,}
	 $\{u^k:=(y^k,\by^k,z^k,\bz^k,x^k,\bx^k)\}$ be the sequence generated by DBA-sGS-ALM. Assume the KKT point set is nonempty. Then, we have the following results.
	
	\begin{enumerate}
		\item[(a)] The sequence $\{(y^k,\by^k,z^k,\bz^k)\}$ converges to an optimal solution of the dual problem (D) of \eqref{eq-DBA}
		and the sequence $\{(x^k,\bx^k)\}$ converges to an optimal solution of the  primal problem (P) of \eqref{eq-DBA}.
		
		\item[(b)] \blue{Suppose that the sequence $\{u^k\}$ converges to a KKT point $\hat{u}:=(\hat{y},\hat{\by},\hat{z},\hat{\bz},\hat{x},\hat{\bx})$ and the KKT mapping $\mathcal{R}$ is metrically subregular at $(\hat{u},0)\in {\rm gph}\mathcal{R}$. If a relative error criterion is employed instead:
		\[
			\|\Delta^k\|\leq\eta_k\|u^k-u^{k+1}\|,
		\]
		 then the sequence $\{u^k\}$ is locally linearly convergent to $\hat{u}$.}
	\end{enumerate}
	
\end{theorem}

%
%
%

\subsection{The efficient computation of $\bar{y}^{k+1}_{\rm tmp}$  and $\by^{k+1}$}
{As $A$ is a single block matrix, the computation of $y^{k+1}$ and $y_{\rm tmp}^{k+1}$ involve solving the linear system of the form $(AA^*+J)y=h$, which is usually not an issue. For example, if computing the Cholesky factorization of $AA^*$ can be done cheaply, we can simply choose $J=0$. If factorizing $AA^*$ is expensive, we can choose $J=\lambda_{\max}(AA^*)I-AA^*$ to efficiently  compute $y$ at the expense of perturbing the augmented Lagrangian function
by the non-trivial proximal term $\frac{\sigma}{2}\norm{y-y^k}_J^2$.

On the other hand, since $B$ and $\bar B$ are multiblock matrices, to implement DBA-sGS-ALM or DBA-sGS-ADMM efficiently, it is necessary to carefully address the efficient computation of $\by^{k+1}_{\rm tmp}$ and $\by^{k+1}$.}
To compute the solution $\by^{k+1}_{\rm tmp}$ in \eqref{eq-step2}, 
we need to solve a generally very large $\bm\times \bm$ system of linear equations of the form: 
\begin{eqnarray}
M \by= h,
\label{eq-M2}
\end{eqnarray}
where the $\bm\times \bm$ matrix $M$ given in \eqref{eq-M} has the following structure:
\begin{eqnarray}
M = \left[ \begin{array}{c} B_1 \\ \vdots \\ B_N \end{array} \right]
\left[ \begin{array}{c} B_1 \\ \vdots \\ B_N \end{array} \right]^* + 
\left[ \begin{array}{ccc} \bB_1\bB_1^* & &  \\ & \ddots &  \\  & & \bB_N\bB_N^*\end{array} \right] + \bJ.
\label{eq-MM}
\end{eqnarray}

Here we design several numerical strategies which are crucial for cutting down the computational cost
of solving \eqref{eq-M2}. 

\begin{description}
	\item[(a)] Suppose we choose $\bJ = \diag{\bJ_1,\ldots,\bJ_N}$ where each $\bJ_i:\cX_i\to\cX_i$ is a 
	self-adjoint positive semidefinite linear operator.
	To compute the inverse of $M$ in \eqref{eq-MM}, we may make use of the Sherman-Morrison-Woodbury (SMW) formula which we shall describe next. 
	Let $\bD = \diag{\bD_1,\dots,\bD_N}$, where $\bD_i = \bB_i\bB_i^*+\bJ_i$, $i=1,\ldots,N.$
	Assume that $\bD_i$ is invertible for all $i=1,\ldots,N$, then 
	we have that
	\begin{eqnarray*}
		&M^{-1} h= (\bD+ BB^*)^{-1} h
		= \bD^{-1} h  - \bD^{-1} B G^{-1}B^*\bD^{-1} h,
		&
	\end{eqnarray*}
	where $G:\cX\to \cX$ is given by
	\begin{eqnarray}
	G := I_n +\sum_{i=1}^N B_i^*\bD_i^{-1}B_i . 
	\label{eq-G}
	\end{eqnarray}
	Thus to solve \eqref{eq-M2}, we only {need to} solve a linear system of the form 
	$Gx =g := B^*\bD^{-1}h.$
	
	Observe that to compute $G$,  the component matrices $B_i^*\bD_i^{-1}B_i $, $i=1,\ldots,N$, can be computed in
	parallel, and operations such as $\bD^{-1}_i h_i$, $i=1,\ldots,N$, can also be done
	in parallel. Similarly,  operations such as evaluating $Bx$ and $B^*\by$ for given $x\in \cX$ and $\by\in\bcY$ can 
	also be done in parallel.
	
	(i) Generally, one would expect the $n\times n$ matrix $B_i^*\bD_i^{-1}B_i$ to be dense even if $B_i$ is sparse, 
	unless $\bD_i$ has a nice structure such as being a diagonal matrix. If the memory
	required to store the matrix $G$ and its Cholesky 
	factor is prohibitive, then one may choose $\bJ_i$ such that 
	$\bD_i$ is a diagonal matrix. For example, we may choose
	$$
	\bJ_i = \lam_{\max}(\bB_i\bB_i^*) I_{m_i} - \bB_i\bB_i^*, \quad i=1,\ldots,N.
	$$
	Then $\bD_i = \lam_{\max}(\bB_i\bB_i^*)  I_{m_i}$. 
	{In which case, $G = I_n+\sum_{i=1}^N \lam_{\max}(\bB_i\bB_i^*) ^{-1} B_i^*B_i$, and it  
	is more likely to be sparse. As an example, consider the case where all $B_i$ are identical and $B_i^*B_i$ is sparse. In this scenario, $G$ will be sparse as well. However, if $\bar J_i=0$ and $(\bar B_i\bar B^*_i)^{-1}$ is dense, then $G$ may become dense even if $B_i^*B_i$ is sparse.}
	{This technique of attempting
		to obtain a sparser $G$ also explains the rationale for introducing the 
		linear operator $\bJ$ in Step 2 of DBA-sGS-ADMM and Step 1 of DBA-sGS-ALM.}
	
	(ii) In general, we can expect the symmetric positive definite matrix $G$ in \eqref{eq-G} to be 
	well conditioned, especially if for each $i$, we choose $\bJ_i = \alpha_i I_{m_i}$ for a {sufficiently large 
	positive} scalar
	$\alp_i$. 
	In this case, a preconditioned conjugate gradient (PCG) method  applied 
	to solve $G x = g$ is expected to converge in a small number of steps. 
	In addition, for a given $x$, the matrix-vector product $Gx$ can also be computed very efficiently 
	since the component vectors, $\bB_i^* \bD_i^{-1} \bB_i x$, $i=1,\ldots,N$, can be computed in parallel.
	

	\item[(b)]
	It is also possible to solve \eqref{eq-M2} in a parallel fashion if we choose $\bJ$ appropriately. 
	For example, if we choose
	\begin{eqnarray}
	\bJ = \diag{B_i B_i^* + \mbox{$\sum_{j=1,j\not= i}^N$} \norm{B_i B_j^*}_2 \,I_{m_i} \mid i=1,\ldots,N} - BB^* \;\succeq\; 0,
	\label{eq-bJ}
	\end{eqnarray}
	then 
	\begin{eqnarray} 
	M = \diag{\bE,\ldots,\bE_N},
	\end{eqnarray}
	where $\bE_i = \bB_i\bB_i^* + B_iB_i^*+ \mbox{$\sum_{j=1,j\not= i}^N$} \norm{B_i B_j^*}_2 I_{m_i},$  $i=1,\ldots,N$. 
	Then the solution $\by$ in \eqref{eq-M2} can be computed in parallel by solving
	\begin{eqnarray*}
		\bE_i \by_i = h_i, \quad i=1,\ldots,N.
	\end{eqnarray*}
	{We should mention that in the literature, for example, \cite{LST-PBA} and references therein}, a standard way to decompose $M$ into 
	a block diagonal form is to choose $\bJ$ to be:
	\begin{eqnarray*}
		\bJ^{\rm std} = (N+1)\diag{B_1B_1^*,\ldots,B_NB_N^*} - BB^* \succeq 0,
	\end{eqnarray*}
	where the positive semidefiniteness of $\bJ^{\rm std}$ can be shown by using the inequality 
	that $\norm{\sum_{i=1}^N \bx_i }^2 \leq N\sum_{i=1}^N\norm{\bx_i}^2$ for any $\bx\in \bcX$. 
	However, such a choice is usually more conservative than the one in \eqref{eq-bJ}, 
	especially when $\norm{B_iB_j^*}_2 \ll 1$ for many pairs in $\{(i,j) \mid 1\leq i<j\leq N\}$.
	
	\item[(c)] We should also take special consideration when the optimization problem \eqref{eq-DBA}
	has the property that  $B_1 = B_2 = \cdots =  B_N$. 
	In this case 
	$$
	G = I_n +B_1^*  \big(\mbox{$\sum_{i=1}^N$}  \bD_i^{-1} \big) B_1.
	$$
	Thus if $m_1$ is not too large, one can precompute the 
	$m_1\times m_1$ matrix $\sum_{i=1}^N\bD_i^{-1}$  once at the beginning of the algorithm
	so that the matrix-vector product $Gx$ can be evaluated very efficiently for any 
	given $x\in \cX.$ 
	If in addition, we also have that $\bB_1= \cdots = \bB_N$ such as in 
	the case of the DNN relaxation 
	\eqref{eq-UFL-DNN} of the UFL problem, then we get 
	$$
	G = I_n + N B_1^* \bD_1^{-1} B_1. 
	$$
	In this case, obviously there is tremendous saving in  the computation of the matrix-vector product $Gx$
	for any given $x$.

\end{description}

\section{Incorporating a semismooth Newton method in DBA-sSG-ADMM and DBA-sGS-ALM
	for solving (D)}

\begin{description}
	
	\item[(a)]
	Recall that in computing $(z^{k+1},y^{k+1})$ in \eqref{eq-yz} for Algorithm DBA-sGS-ADMM, we have added a
	proximal term based on the sGS linear operator ${\bf sGS}(\cQ_{z,y})$ to 
	decouple the computation of $z^{k+1}$ and $y^{k+1}$. 
	We should note that while the computation of the individual $z^{k+1}$ and $y^{k+1}$ 
	has been made easier, the negative effect of adding a possibly large proximal term is that the 
	overall algorithm may take more iterations to converge to a required level of accuracy. 
	Thus in practice, a judicious choice must always be made to balance the
	two competing factors. In particular, if an efficient solver is available to 
	compute $(z^{k+1},y^{k+1})$ simultaneously in \eqref{eq-yz} 
	without the need to add a proximal term, i.e., $J=0$ and $\cS=0$, then 
	one should always adopt this option.
	In this case, we have that 
	\begin{eqnarray}
	(z^{k+1},y^{k+1}) =
	\argmin\limits_{z,y}  \left\{
	\begin{array}{l}
	-\inprod{b}{y} + \delta_\cK^*(-z) -{\inprod{\delta^{k}}{(z,y)}} \\[5pt]
	+
	\frac{\sigma}{2}\norm{A^*y + z + B^*\by^k+v^k-c^k}^2 
	\end{array}
	\Bigm\vert\; \begin{array}{l}  y\in\cY, \\[3pt]  z \in \cX \end{array}
	\right\}.
	\quad \label{eq-ssn1}
	\end{eqnarray}

	{In Step 1 of Algorithm DBA-sGS-ALM, to compute $(\bar z^{k+1},z^{k+1},y^{k+1},\bar y^{k+1})$ efficiently in sequential manner, we regard $(\bar z,z,y)$ and $\bar y$ in \eqref{eq-ssn1} as $x_1$ and $x_2$ in the sGS decomposition theorem, respectively.} By choosing $J=0$, and $\cS = {\bf sGS}(\cQ_{\bz z y,\by})$, the sGS linear operator associated with the following positive semidefinite linear operator,
	\begin{eqnarray*}
		\cQ_{\bz z y, \by} = \left[ 
		\begin{tabular}{ccc|c} 
			$I_{\bn}$   & 0   & 0   & $\bB^*$
			\\ 
			0 & $I_{n}$  &  $A^*$   & $B^*$ 
			\\ 
			0 & $A$    & $AA^* $    &$AB^*$
			\\ \hline
			\blue{\small{$\bB$}} & $B$   &$BA^*$ &\small{$BB^* + \blue{\bB\bB^*} + \bar{J}$}
		\end{tabular}
		\right],
	\end{eqnarray*}
	then we can compute $(\bz^{k+1},z^{k+1},y^{k+1},\by^{k+1})$
	in a symmetric Gauss-Seidel fashion as follows.
	Let  $c^k = c-\sig^{-1} x^k$ and 
	$\bc^k = \bc - \sig^{-1}\bx^k.$ Compute
	{\small
	\begin{eqnarray}
	&&\by^{k+1}_{\rm tmp} \; =\; \argmin\limits_{\by}  \left\{
	\begin{array}{l}
	-\inprod{\bb}{\by}
	+\frac{\sigma}{2}\norm{B^*\by  + A^* y^{k} + z^k - c^k}^2 
	+\frac{\sigma}{2}\norm{\bB^*\by + \bz^k - \bc^k}^2 
	\\[5pt]
	+\frac{\sig}{2}\norm{\by-\by^k}_{\bJ}^2 -\inprod{\bar{\delta}^{k}_{\rm tmp}}{\by}
	\end{array}
	\right\}
	\nn \\[8pt]
	& &\left\{  \begin{array}{l}
	(z^{k+1},y^{k+1}) \; = \;
	\argmin\limits_{z,y}  \left\{ 
	\begin{array}{l} -\inprod{b}{y} + \delta_\cK^*(-z)+
	\\[5pt]
	\frac{\sigma}{2}\norm{z+A^*y + B^*\by_{\rm tmp}^{k+1}-c^{k}}^2 
	-{\inprod{\delta^{k}}{(z,y)}}
	\end{array}
	\Bigm\vert \;
	\begin{array}{l} 
	y\in \cY, \\ z \in \cX
	\end{array} 
	\right\}
	\\[18pt]
	\bz^{k+1} \;=\; 
	\argmin\limits_{\bz}  \big\{
	\frac{\sigma}{2}\norm{\bz+ \bB^*\by_{\rm tmp}^{k+1}-\bc^k}^2  +\delta_{\bcK}^*(-\bz) \mid \bz \in \bcX
	\big\}
	\end{array}\right.
	\qquad \label{eq-ssn2} 
	\\[8pt]
	&& \by^{k+1}\; = \; \argmin\limits_{\by}  \left\{
	\begin{array}{l}
	-\inprod{\bb}{\by}
	+\frac{\sigma}{2}\norm{B^*\by+A^*y^{k+1}+ z^{k+1}-c^k}^2
	+\frac{\sigma}{2}\norm{\bB^*\by+ \bz^{k+1}-\bc^k}^2 
	\\[5pt]
	+\frac{\sig}{2}\norm{\by-\by^k}_{\bJ}^2 - \inprod{\bar{\delta}^{k}}{\by}
	\end{array}
	\right\}.
	\nn
	\end{eqnarray}
	}
	
	In both \eqref{eq-ssn1} and \eqref{eq-ssn2}, we need to solve a problem of the form:
	{\small
	\begin{eqnarray}
	(z^{k+1},y^{k+1}) \; = \; 
	\argmin\limits_{z,y}  \Big\{ -\inprod{b}{y} + \delta_\cK^*(-z) + 
	\frac{\sigma}{2}\norm{z+A^*y - \widehat{c}^k}^2 -{\inprod{\delta^{k}}{(z,y)}}
	\Bigm\vert \begin{array}{l} 
	y\in \cY, \\ z \in \cX
	\end{array} 
	\Big\}
	\quad
	\label{ALMSNN-subp-zy}
	\end{eqnarray}
	}
	for some given $\widehat{c}^k$. 
	By making use of the projection onto $\cK$, we can first compute $y^{k+1}$ from the 
	following minimization problem, and obtain $z^{k+1}$ once $y^{k+1}$ is computed as follows:
	\begin{eqnarray}
		&& \left\{ \begin{array}{rcl}
			y^{k+1} & \approx & 
			\argmin\limits_{y}  \left\{
			\begin{array}{l}
				-\inprod{b}{y} +
				\sigma M_{\delta_{\cK/\sigma}^*}(A^*y -\widehat{c}^k)-\inprod{\delta_y^k}{y}
			\end{array}
			\bigm\vert\; y\in\cY
			\right\}
			\\[5pt]
			z^{k+1} &=& \frac{1}{\sig}\Pi_\cK(\sig(A^*y^{k+1}-\hat c^k))-(A^*y^{k+1}-\hat c^k).
		\end{array} \right.
		\label{ALMSNN-subp2-zy}
	\end{eqnarray}
	\blue{The detail of deriving \eqref{ALMSNN-subp2-zy} can be found in Appendix 2.}
		In the interest of keeping the paper to be within a
		reasonable length, we do not delve into the details
		on how one can compute $y^{k+1}$ via a semismooth Newton method, {interested readers can refer to, for example, \cite{LSTc} for a detailed explanation on how to use a semismooth Newton method.} 
	
	\item[(b)]
	In the second step in DBA-sGS-ADMM, if we further remove the proximal term with respect to $(v,\bar v,\bar y)$, i.e. $\bar \cJ=0,\cT=0$, we compute
	\begin{eqnarray}
	\hspace{-5mm}
	(v^{k+1},\bar v^{k+1},\bar y^{k+1}) =
	\argmin\limits_{v,\bv,\by} \left\{
	\begin{array}{l}
	\theta^*(-v)+\bar\theta^*(-\bar v)-\inprod{\bar b}{\bar y} + 
	\frac{\sigma}{2}\norm{ B^*\by+z^{k+1}+\bar v-\bar c^k}^2 \\[3pt]
	+\frac{\sigma}{2}\norm{A^*y^{k+1} + B^*\by+z^{k+1}+v-c^k}^2
	-\inprod{\bar{\delta}^k}{\bar{y}}
	\end{array}
	\right\}.
	\quad \label{eq-ssn3}
	\end{eqnarray}
	Let $R^{5,k}=A^*y^{k+1}+z^{k+1}-\hat c^k$, $R^{6,k}=\bar z^{k+1}-\hat{\bar c}^k$, then solution of \eqref{eq-ssn3} can be obtained as follows:
	\begin{eqnarray}
		&& \left\{ \begin{array}{rcl}
			\by^{k+1} & \approx & 
			\argmin\limits_{\by\in\bar\cY} \left\{
			\begin{array}{l}
				-\inprod{\bar b}{\bar y} +
				\sigma M_{\theta^*/\sigma}(R^{5,k}+B^*\bar y)+\sigma M_{\bar\theta^*/\sigma}(R^{6,k}+\bar B^*\bar y)-\inprod{\bar\delta^k}{\bar y}
			\end{array}
			\right\}
			\\[5pt]
			\bar v^{k+1} &=& -\mathop{\rm Prox}_{\bar\theta^*/\sigma}(R^{6,k}+\bar B^*\bar y^{k+1})
			\\[5pt]
			v^{k+1} &=& -\mathop{\rm Prox}_{\theta^*/\sigma}(R^{5,k}+B^*\bar y^{k+1}).
		\end{array} \right.
		\label{ALMSNN-subp2-yv}
	\end{eqnarray}
	\blue{The detail of deriving \eqref{ALMSNN-subp2-yv} can be found in Appendix 2.}
	The scale of the subproblem solving for $\bar{y}^{k+1}$ can be large when $N$ increases. On the other hand, the update scheme of $\bar v^{k+1}$ shows that this subproblem can be solved in a parallel manner. We leave the efficient computation of these two subproblems as a future work.
	
\end{description}


\section{Numerical experiments}


In this section, we conduct numerical experiments to test the efficiency of our algorithm against state-of-the-art solvers. We use data sets found in the literature (such as two-stage stochastic programming problems and uncapacitated facility location problems) as well as synthetic data. All the experiments are run in \texttt{MATLAB} R2017b on a computer with 1.4 GHz Quad-Core Intel Core i5 Processor (this processor has 4 cores and 8 threads) and 8 GB memory, equipped with macOS.

\subsection{Stopping conditions}
Based on the optimality conditions in \eqref{eq-KKT}, we measure the accuracy of a computed solution by the following relative KKT residue:
\[
\eta = \max\{\eta_P,\eta_D,0.2\eta_{\cK},\blue{0.2\eta_{\theta}},\eta_{\bar{P}},\eta_{\bar{D}},0.2\eta_{\bcK},\blue{0.2\eta_{\bar\theta}}\},
\]
\noindent 
where
\begin{eqnarray*}
	\begin{array}{cccc}
		\eta_P = \frac{\|Ax-b\|}{1+\|b\|}, &
		\eta_D = \frac{\|A^*y + \bB^* \by + z + v - c\|}{1+\|c\|}, &
		\eta_{\cK} = \frac{\|x-\Pi_{\cK}(x-z)\|}{1+\|x\|+\|z\|}, &\blue{{\eta_\theta=\frac{\|x-\mathop{\rm{Prox}}_\theta(x-v)\|}{1+\|x\|+\|v\|}}}, \\[5pt]
		\eta_{\bar{P}} = \frac{\|Bx + \bB \bx - \bb\|}{1+\|\bar b\|}, &
		\eta_{\bar{D}} = \frac{\|\bB^* \by + \bz + \bv - \bc\|}{1+\|\bc\|}, &
		\eta_{\bcK} = \frac{\|\bx-\Pi_{\bcK}(\bx - \bz)\|}{1+\|\bx\|+\|\bz\|},& \blue{\eta_{\bar\theta}=\frac{\|\bx-\mathop{\rm{Prox}}_{\bar\theta}(\bx-\bv)\|}{1+\|\bx\|+\|\bv\|}}.	
	\end{array}
\end{eqnarray*}
\noindent 
In addition, we also compute the duality gap by:
\blue{
\begin{eqnarray*}
\eta_{\rm gap} = \frac{{\rm obj}_{\rm P} - {\rm obj}_{\rm D}}{1+\mid{\rm obj}_{\rm P}\mid+\mid{\rm obj}_{\rm D}\mid},
\end{eqnarray*}}
where 
\begin{eqnarray*}
	{\rm obj}_{\rm P} &:=& \theta(x) + \inprod{c}{x} 
	+ \btheta(\bx)+\inprod{\bc}{\bx}, \\
	{\rm obj}_{\rm D} &:=& -\theta^*(-v) - \delta_\cK^*(-z) 
	+\inprod{b}{y} - \btheta^*(-\bar{v}) - \delta_{\bcK}^*(-\bz) +\inprod{\bb}{\by},
\end{eqnarray*}
are the primal and dual objective values, respectively.
We terminate our algorithm when $\eta \le 10^{-5}$ and $\eta_{\rm gap}\le 10^{-4}$.

\subsection{2-stage stochastic programming problems}

In this section, we aim to solve two-stage stochastic programming problems on some real data sets. We solve linear problems in the first subsection and quadratic problems in the second subsection. Since there are not much quadratic stochastic programming test data available, we will follow the idea mentioned in \cite{LauWomersley} to extend stochastic linear programming problems to its quadratic counterparts by adding a small quadratic term in the objective. 
Specifically, we choose to set
$
\theta(x) = \frac{0.1}{2}\norm{x}^2, \quad \btheta(\bx) = \frac{0.1}{2}\norm{\bx}^2.
$
\blue{$\cK$ and $\bar\cK$ in these instances are box constraints.}


The detailed description of each data set used in this section can be found in the Appendix. 
In table \ref{table:SPdatasets}, we only list down the name and source of the data sets.

\begin{table}
	\begin{tabular}{|c|c|}
		\hline 
		Name & Source \\ \hline 
		\begin{tabular}{c}
			\texttt{assets} \\
			\texttt{env} \\
			\texttt{phone} \\
			\texttt{AIRL}\\
			\texttt{4node}
		\end{tabular} & 
		\url{http://www4.uwsp.edu/math/afelt/slptestset/download.html} \\ \hline 
		\begin{tabular}{c}
			\texttt{pltexp}\\ \texttt{storm}
		\end{tabular} & 
		\url{http://users.iems.northwestern.edu/~jrbirge/html/dholmes/post.html} \\ \hline 
		\texttt{gbd} & \url{http://pages.cs.wisc.edu/~swright/stochastic/sampling/} \\ \hline 
	\end{tabular}
	\caption{Source of two-stage stochastic programming data sets.} 
	\label{table:SPdatasets}
\end{table}


We compare our proposed DBA-sGS-ADMM with Gurobi under its default setting 
\blue{except for modifying the accuracy tolerance parameters: {{\tt FeasibilityTol, OptimalityTol, BarConvTol}\footnote{\url{https://www.gurobi.com/documentation/9.5/refman/parameters.html.}}, which are absolute primal feasibility, absolute dual feasibility and relative gap between primal and dual values respectively. More specifically, we set {\tt FeasibilityTol}=$\left(1+\max\{\|b\|,\|\bar b\|\}\right)\times10^{-5}$, {\tt OptimalityTol}=$\left(1+\max\{\|c\|,\|\bar c\|\}\right)\times10^{-5}$, and {\tt BarConvTol}=$10^{-4}$.} In this setting, Gurobi adopts concurrent optimizers\footnote{\url{https://www.gurobi.com/documentation/9.5/refman/concurrent\_optimizer.html}} for linear programming, which runs multiple solvers (the primal/dual simplex method and the barrier method) simultaneously, and chooses the one that finishes first. 
Concurrent strategy requires multiple (at least 3) threads, and the first three threads are assigned to primal simplex, dual simplex and barrier respectively. The extra threads are then assigned to the barrier method. More threads will not make the simplex method any faster. However, the matrix factorization operation in the barrier method is faster with more threads. We use 8 threads for Gurobi in order to adopt the concurrent strategy for linear programming. Therefore, when solving linear programming, our method is almost compared with the best of the simplex method and the barrier method, meanwhile the barrier method is accelerated by multiple threads. Even though this comparison is a bit unfair to our method, the results later show that the performance of our method is comparable to Gurobi, and even better than Gurobi for some large-scale instances.} 

For quadratic programming, Gurobi uses the barrier method by default. In this case, it does improve the efficiency with more threads, but when the number of threads increases, the speed is not necessarily faster because the cost of data transmission will also increase. Meanwhile, 
	the computational efficiency of our DBA-sGS-ADMM/ALM can be further improved by incorporating parallel computation. Hence in quadratic programming, we only use a single thread in Gurobi for a fair comparison.
	
	We incorporate a semismooth Newton method in DBA-sGS-ADMM for some examples for which the size of the matrix $A$ is small while the number of scenarios is relatively large. More specifically, we adopt a semismooth Newton method described in Section 6 when $m_0\leq10,n_0\leq20$ and $N\geq100$ in our numerical experiments.

We adopt different strategies to solve the linear system \eqref{eq-M2}  efficiently for different cases. For some simple cases, such as when the {component blocks of $B$ or $\bar B$
		are identical,} or the data size is not very large (more specifically, when $\sum_{i=1}^N m_i<5000$), we directly solve the linear system \eqref{eq-M2} by utilizing Cholesky factorization. When the data size is relatively large (specifically $\sum_{i=1}^N m_i\geq5000$), we make use of SMW formula described in Section 5.4 to solve the linear system. In our numerical experiment, when $\sum_{i=1}^N m_i<50000$, we use the SMW formula without the proximal term $\bJ$, i.e. $\bar J=0$. When the data scale is very large, more specifically $\sum_{i=1}^N m_i\geq50000$, we use the SMW formula by setting $\bar J_i=\lambda_{\max}(\bar B_i\bar B_i^*)I_{m_i}-\bar B_i\bar B_i^*$, $i=1,...,N$ to make the matrix 
	$G$ in \eqref{eq-G} more likely to be sparse.

For stochastic linear programming problems, we also compare our algorithm with the well-known progressive hedging algorithm (PHA), which is applied to the primal problem (P). We should emphasize that just like the PHA, the subproblems associated with the second stage variables $\bar y,\bar z,\bar v$ in our algorithm can also be solved in a parallel manner, as mentioned in section 5 and section 6(b). However, in our case, the cost of each parallel computation is already very cheap. If we implement the distributed version using \texttt{MATLAB}'s parallel for-loop, the overhead cost will become dominant and affect the performance adversely. Hence, we choose to implement the algorithm serially in this paper. For a fair comparison, it is reasonable to compare our non-parallel version with the non-parallel version of PHA. 
Recent work in \cite{SXZ} has shown that one can choose the step length $\tau$ in PHA
to be in the range $(0,(1+\sqrt{5})/2)$. In our numerical experiments, we set
$\tau=1.618$ for both 
DBA-sGS-ADMM and PHA. 

\subsubsection{2-stage stochastic linear programming problems}
\def\objP{$\mbox{obj}_{\rm P}$}

\blue{
\begin{center}
	\begin{tiny}       
		\setlength\LTcapwidth{0.9\linewidth}
		\setlength{\tabcolsep}{0.9pt}
		\begin{longtable}{|| c | c | c | c || c | c | c || c | c | c || c | c | c |}
			\hline
			\multicolumn{4}{||c||}{} & \multicolumn{3}{|c||}{DBA-sGS-ADMM} & \multicolumn{3}{|c||}{Gurobi}  & \multicolumn{3}{|c|}{PHA}\\ \hline
			\multicolumn{1}{||c|}{Data} & \multicolumn{1}{|c|}{$m_0 \mid m_i$} & \multicolumn{1}{|c|}{$n_0 \mid n_i$}
			& \multicolumn{1}{|c||}{$N$}  & \multicolumn{1}{|c|}{Iter}&\multicolumn{1}{|c|}{Time(s)}  & \multicolumn{1}{|c||}{\objP}  &\multicolumn{1}{|c|}{Iter} &\multicolumn{1}{|c|}{Time(s)} & \multicolumn{1}{|c||}{\objP}  &\multicolumn{1}{|c|}{Iter} &\multicolumn{1}{|c|}{Time(s)} & \multicolumn{1}{|c|}{\objP} 
			\\ \hline
			\endhead
			
			phone\textunderscore 1 &   1 $\mid$   23 &   9 $\mid$   93 &   1 & 2301 &0.33 &3.69e+01 & 25 &0.03 &3.69e+01 &153 &0.37 &3.69e+01
			\\[3pt]\hline
			AIRL\textunderscore first &  2 $\mid$    6 &   6 $\mid$   12 &  25 &12651 &1.73  &2.49e+05 & 146 &0.01 &2.49e+05 &312 &11.99 &2.49e+05
			\\[3pt]\hline
			AIRL\textunderscore second &   2 $\mid$    6 &   6 $\mid$   12 &  25 &14751 &2.02  &2.70e+05 & 157 &0.01 &2.70e+05 &242 &9.28 &2.70e+05
			\\[3pt]\hline
			4node\textunderscore 256 &  14 $\mid$   74 &  60 $\mid$  198 & 256 & 637 &2.4  &4.25e+02 & 25078 &2.72 &4.25e+02 &504 &492.39 &4.25e+02
			\\[3pt]\hline
			4node\textunderscore 512 &  14 $\mid$   74 &  60 $\mid$  198 & 512 & 2681 &7.54  &4.30e+02 & 49908 &8.37 &4.30e+02 &813 &1617.74 &4.30e+02
			\\[3pt]\hline
			4node\textunderscore 1024 &  14 $\mid$   74 &  60 $\mid$  198 &1024 & 1450 &7.89  &4.34e+02 & 119424 &11.10 &4.34e+02 &1022 &3945.46 &4.34e+02
			\\[3pt]\hline
			4node\textunderscore 2048 &  14 $\mid$   74 &  60 $\mid$  198 &2048 & 1287 &13.88  &4.42e+02 & 259486 &27.42 &4.42e+02 &-&-&-
			\\[3pt]\hline
			4node\textunderscore 4096 &  14 $\mid$   74 &  60 $\mid$  198 &4096 & 1628 &36.64  &4.47e+02 & 464119 &66.96 &4.47e+02 &-&-&-
			\\[3pt]\hline
			4node\textunderscore 8192 &  14 $\mid$   74 &  60 $\mid$  198 &8192 & 767 &41.14  &4.47e+02 & 141 &239.78 &4.47e+02 &-&-&-
			\\[3pt]\hline
			4node\textunderscore 16384 &  14 $\mid$   74 &  60 $\mid$  198 &16384 & 782 &89.88  &4.47e+02 & 155 &474.69 &4.47e+02 &-&-&-
			\\[3pt]\hline
			4node\textunderscore 32768 &  14 $\mid$   74 &  60 $\mid$  198 &32768 & 748 &186.84  &4.47e+02 &164 &869.61 &4.47e+02 &-&-&-
			\\[3pt]\hline
			4node\textunderscore base\textunderscore 1024 &  16 $\mid$   74 &  60 $\mid$  198 &1024 & 1881 &10.39  &4.23e+02 & 121 &20.36 &4.23e+02 &-&-&-
			\\[3pt]\hline
			4node\textunderscore base\textunderscore 2048 &  16 $\mid$   74 &  60 $\mid$  198 &2048 & 1511 &19.01  &4.24e+02 & 125 &53.30 &4.25e+02 &-&-&-
			\\[3pt]\hline
			4node\textunderscore base\textunderscore 4096 &  16 $\mid$   74 &  60 $\mid$  198 &4096 & 1290 &33.63  &4.25e+02 & 127 &112.69 &4.25e+02 &-&-&-
			\\[3pt]\hline
			4node\textunderscore base\textunderscore 8192 &  16 $\mid$   74 &  60 $\mid$  198 &8192 & 898 &54.21  &4.25e+02 &158 &326.67 &4.25e+02 &-&-&-
			\\[3pt]\hline
			4node\textunderscore base\textunderscore 16384 &  16 $\mid$   74 &  60 $\mid$  198 &16384 & 968 &122.31  &4.25e+02 &117 &522.26 &4.25e+02 &-&-&-
			\\[3pt]\hline
			4node\textunderscore base\textunderscore 32768 &  16 $\mid$   74 &  60 $\mid$  198 &32768 & 840 &234.75  &4.25e+02 &152 &1510.31 &4.25e+02 &-&-&-
			\\[3pt]\hline

			\caption{
				Comparison of computational results between DBA-sGS-ADMM, Gurobi and PHA for two-stage stochastic linear programming problems. The results of Gurobi are obtained using \textbf{8 threads}. The objective function value is shown to help us verifying the correctness easily. ``-'' means that the solver runs out of time.
			}
			\label{table:SLP}
		\end{longtable}
	\end{tiny} 
\end{center}
}

From Table \ref{table:SLP}, we can observe that {DBA-sGS-ADMM outperforms Gurobi especially when the number of scenarios $N$ is large.}
{For example, for the instance {\tt 4node\_32768} with about $3.27\times 10^4$ scenarios, we are faster by  at least 4 times.}
We also see that PHA is the least efficient, and it can only solve moderate scale problems, while DBA-sGS-ADMM can solve problem with up to $3.27\times10^4$ scenarios and medium-size dimensional decision variables. The inefficiency of the PHA could be due to the fact that at each 
iteration, it needs to solve $N$ smaller LP subproblems whose total computational
cost are  rather high. {{We have to mention that under the framework of PHA, the cost of computing the projections onto the feasible sets of the scenarios subproblems to check the primal feasibility 
		condition can be quite expensive.}
	Therefore, in practice, one only checks the feasibility of the nonanticipativity constraint and the relative error of iteration points or function values. Therefore, {another advantage of} our algorithm compared {to} PHA is that the stopping criterion (KKT residue) of DBA-sGS-ADMM is 
	theoretically more 
	rigorous than that of the PHA.}





\subsubsection{2-stage stochastic quadratic programming problems}

\blue{
\begin{center}
	\setlength\LTcapwidth{0.9\linewidth}
	\setlength{\tabcolsep}{2pt}
	\tiny
	\begin{longtable}{|| c | c | c | c ||  c | c | c || c | c | c |}
		\hline
		\multicolumn{4}{||c||}{} & \multicolumn{3}{|c||}{DBA-sGS-ADMM} & \multicolumn{3}{|c|}{Gurobi} \\ \hline
		\multicolumn{1}{||c|}{Data} & \multicolumn{1}{|c|}{$m_0 \mid m_i$} & \multicolumn{1}{|c|}{$n_0 \mid n_i$}
		& \multicolumn{1}{|c||}{$N$}  & \multicolumn{1}{|c|}{Iter}&\multicolumn{1}{|c|}{Time(s)}  & \multicolumn{1}{|c||}{\objP} 	
		&\multicolumn{1}{|c|}{Iter} &\multicolumn{1}{|c|}{Time(s)} & \multicolumn{1}{|c|}{\objP}  \\ \hline
		
		\endhead
		
		assets\textunderscore small &   5 $\mid$    5 &  13 $\mid$   13 & 100 &577 &0.47  &1.86e+05 & 21 &0.43 &1.86e+05
		\\[3pt]\hline
		assets\textunderscore large &   5 $\mid$    5 &  13 $\mid$   13 &37500 &6801 &395.38  &4.92e+06 & 27 &14.12 &4.92e+06
		\\[3pt]\hline
		env\textunderscore aggr &  32 $\mid$   48 &  69 $\mid$   85 &   5 & 501 &0.49  &6.38e+04 & 23 &0.20 &6.38e+04
		\\[3pt]\hline
		env\textunderscore loose &  32 $\mid$   48 &  69 $\mid$   85 &   5 & 255 &0.16 
		&8.07e+04 & 21 &0.49 &8.07e+04
		\\[3pt]\hline
		env\textunderscore first &  32 $\mid$   48 &  69 $\mid$   85 &   5 & 423 &0.17 
		&6.93e+04 & 21 &0.41 &6.93e+04
		\\[3pt]\hline
		env\textunderscore imp &  32 $\mid$   48 &  69 $\mid$   85 &  15 &1022 &0.50  &1.22e+05 & 25 &0.20 &1.22e+05
		\\[3pt]\hline
		env\textunderscore 1200 &  32 $\mid$   48 &  69 $\mid$   85 &1200 &2872 &14.38  &6.68e+06 & 36 &3.89 &6.68e+06
		\\[3pt]\hline
		env\textunderscore 1875 &  32 $\mid$   48 &  69 $\mid$   85 &1875 &4800 &21.62  &1.05e+07 & 37 &6.02 &1.05e+07
		\\[3pt]\hline
		env\textunderscore 3780 &  32 $\mid$   48 &  69 $\mid$   85 &3780 &5099 &44.45  &2.11e+07 & 43 &8.98 &2.11e+07
		\\[3pt]\hline
		env\textunderscore 5292 &  32 $\mid$   48 &  69 $\mid$   85 &5292 &5361 &76.59  &2.95e+07 & 42 &11.55 &2.95e+07
		\\[3pt]\hline
		env\textunderscore lrge &  32 $\mid$   48 &  69 $\mid$   85 &8232 &4611 &125.11  &4.58e+07 & 39 &20.17 &4.58e+07
		\\[3pt]\hline
		env\textunderscore xlrge &  32 $\mid$   48 &  69 $\mid$   85 &32928 &4553 &1470.55  &1.85e+08 & 38 &121.71 &1.85e+08
		\\[3pt]\hline
		pltexpA2\textunderscore 6 &  62 $\mid$  104 & 188 $\mid$  272 &   6 &20000 &5.20 &6.50e+04 & 25 &0.48 &6.50e+04
		\\[3pt]\hline
		pltexpA2\textunderscore 16 &  62 $\mid$  104 & 188 $\mid$  272 &  16 &20000 &7.91  &1.54e+05 & 24 &0.58 &1.54e+05
		\\[3pt]\hline
		stormG2\textunderscore 8 &  58 $\mid$  528 & 179 $\mid$ 1377 &   8 &14401 &28.59  &1.83e+07 & 63 &0.96 &1.83e+07
		\\[3pt]\hline
		stormG2\textunderscore 27 &  58 $\mid$  528 & 179 $\mid$ 1377 &  27 &8101 &23.88  &2.11e+07 & 66 &5.12 &2.11e+07
		\\[3pt]\hline
		stormG2\textunderscore 125 &  58 $\mid$  528 & 179 $\mid$ 1377 & 125 &15501 &127.37  &3.04e+07 & 62 &49.37 &3.04e+07
		\\[3pt]\hline
		gbd &   4 $\mid$    5 &  21 $\mid$   10 &646425 &1851 &462.80  &7.94e+08 & 21 &118.87 &8.21e+08
		\\[3pt]\hline
		phone &   1 $\mid$   23 &   9 $\mid$   93 &32768 &501 &48.89  &2.74e+05 &34 &56.63 &2.74e+05
		\\[3pt]\hline
		phone\textunderscore 1 &   1 $\mid$   23 &   9 $\mid$   93 &   1 & 157 &0.08  &4.77e+01 & 12 &0.37 &4.77e+01
		\\[3pt]\hline
		AIRL\textunderscore first &   2 $\mid$    6 &   6 $\mid$   12 &  25 &1551 &0.26  &3.61e+06 & 21 &0.33 &3.61e+06
		\\[3pt]\hline
		AIRL\textunderscore second &   2 $\mid$    6 &   6 $\mid$   12 &  25 &1651 &0.25  &3.55e+06 & 23 &0.35 &3.55e+06
		\\[3pt]\hline
		4node\textunderscore 1 &  14 $\mid$   74 &  60 $\mid$  198 & 1 & 1371 &0.25  &8.95e+03 & 24 &0.31 &8.95e+03
		\\[3pt]\hline
		4node\textunderscore 2 &  14 $\mid$   74 &  60 $\mid$  198 & 2 & 1311 &0.25  &9.02e+03 & 24 &0.40 &9.02e+03
		\\[3pt]\hline
		4node\textunderscore 256 &  14 $\mid$   74 &  60 $\mid$  198 & 256 & 881 &1.10 &1.25e+04 & 24 &1.15 &1.25e+04
		\\[3pt]\hline
		4node\textunderscore 512 &  14 $\mid$   74 &  60 $\mid$  198 & 512 & 764 &1.61  &1.39e+04 & 24 &2.05 &1.39e+04
		\\[3pt]\hline
		4node\textunderscore 1024 &  14 $\mid$   74 &  60 $\mid$  198 &1024 & 691 &2.97 &1.64e+04 & 31 &4.41 &1.64e+04
		\\[3pt]\hline
		4node\textunderscore 2048 &  14 $\mid$   74 &  60 $\mid$  198 &2048 & 763 &6.26  &2.15e+04 & 50 &17.31 &2.15e+04
		\\[3pt]\hline
		4node\textunderscore 4096 &  14 $\mid$   74 &  60 $\mid$  198 &4096 & 873 &17.12  &3.14e+04 & 36 &25.57 &3.14e+04
		\\[3pt]\hline
		4node\textunderscore 8192 &  14 $\mid$   74 &  60 $\mid$  198 &8192 & 912 &42.18  &5.05e+04 & 53 &107.71 &5.05e+04
		\\[3pt]\hline
		4node\textunderscore 16384 &  14 $\mid$   74 &  60 $\mid$  198 &16384 & 839 &84.11   &8.64e+04 & 33 &100.94 &8.65e+04
		\\[3pt]\hline
		4node\textunderscore 32768 &  14 $\mid$   74 &  60 $\mid$  198 &32768 & 646 &142.01  &1.55e+05 &29 &164.33 &1.55e+05
		\\[3pt]\hline
		4node\textunderscore base\textunderscore 1024 &  16 $\mid$   74 &  60 $\mid$  198 &1024 & 501 &2.19  &1.60e+04 & 65 &6.69 &1.60e+04
		\\[3pt]\hline
		4node\textunderscore base\textunderscore 2048 &  16 $\mid$   74 &  60 $\mid$  198 &2048 & 501 &4.65  &2.02e+04 & 48 &61.26 &2.02e+04
		\\[3pt]\hline
		4node\textunderscore base\textunderscore 4096 &  16 $\mid$   74 &  60 $\mid$  198 &4096 & 501 &10.86  &2.85e+04 & 30 &21.29 &2.85e+04
		\\[3pt]\hline
		4node\textunderscore base\textunderscore 8192 &  16 $\mid$   74 &  60 $\mid$  198 &8192 & 811 &46.01  &4.51e+04 &48 &45.12 &4.51e+04
		\\[3pt]\hline
		4node\textunderscore base\textunderscore 16384 &  16 $\mid$   74 &  60 $\mid$  198 &16384 & 931 &111.32  &7.82e+04 &188 &701.29 &7.82e+04
		\\[3pt]\hline
		4node\textunderscore base\textunderscore 32768 &  16 $\mid$   74 &  60 $\mid$  198 &32768 &1271 &313.84  &1.44e+05 & 63 &3886.62 &1.44e+05
		\\[3pt]\hline

		\caption{
			Comparison of computational results between DBA-sGS-ADMM and Gurobi for two-stage quadratic stochastic programming problem (with added 0.1$I$ in the quadratic objective). All the results are obtained using a {\textbf{8 thread}}. The objective function value is shown to help us verifying the correctness easily.}
		\label{table:SP}
	\end{longtable}
	
\end{center}
}

In table \ref{table:SP}, we observe that the computational performance of DBA-sGS-ADMM and Gurobi are comparable for many instances, in which
the runtime of DBA-sGS-ADMM is shorter among $21/38\approx55.3\%$ of the tested instances. 
In particular, for a few of the {\tt 4node\_base\_xxxx} instances, DBA-sGS-ADMM can outperform
Gurobi by a large margin. 
On the other hand, 
Gurobi can solve some of the instances such as {\tt env\_xxxx} very efficiently. 
This is reasonable, since we only add a very sparse term of $0.1I$ in the quadratic objective, and
Gurobi can {take full advantage of the sparsity in the data.} Later in section 7.4.1, we will see that
when the matrix in the quadratic objective is dense, we can expect extreme inefficiency 
of the barrier method used by Gurobi, and DBA-sGS-ADMM has substantial computational 
advantage.

\subsection{DNN relaxations of uncapacitated facility location (UFL) problems}

\subsubsection{DNN relaxations of linear UFL problems}

In this subsection, we aim to solve some DNN relaxations of  linear UFL problems. We have two benchmark datasets of different scales, namely the small scale \texttt{Bilde-Krarup} dataset and the large scale \texttt{Koerkel-Ghosh} dataset. Both of these data instances could be downloaded from the UflLib at \url{http://resources.mpi-inf.mpg.de/departments/d1/projects/benchmarks/UflLib/index.html}. These are all linear UFL problem, i.e. $q_{ij}=0\;\forall i=1,\ldots,p, j=1,\ldots,q$ in 
\eqref{eq-UFL}.

We should emphasize that our intention here is not to investigate the tightness
	of the DNN relaxation \eqref{eq-UFL-DNN} of the UFL problem \eqref{eq-UFL} but to treat the DNN relaxation
	as an example of the DBA problem \eqref{eq-DBA} with constraint set
	$\cK$ being a positive semdefinite matrix cone $\S^{n}_+$.



We first solve the DNN relaxations  of linear UFL problems described in  
\eqref{eq-UFL-DNN}. 
The numerical results obtained by DBA-sGS-ADMM will be compared against MOSEK 
for the 
small scale dataset  \texttt{Bilde-Krarup} and the results are
reported in Table \ref{table:SUFL-SDP}.
Note that MOSEK is the current best solver for solving general SDP problems. When it is employed to solve the dual block angular problem \eqref{eq-UFL-DNN}, 
MOSEK treats it as a general SDP problem since it does not exploit the 
block-angular structure in the constraint matrix except for sparsity in the data matrices.

%

{Note that for UFL problems, we can make use of the special structure of $B$ and $\bar B$ in the constraint to further cut down the cost of solving the linear system \eqref{eq-M2}. 
	Specifically, the component blocks of $B$ and $\bar B$ are repeated, i.e. $B_1=B_2=...=B_N$, $\bar B_1=\bar B_2=...=\bar B_N$. Furthermore, $\bar B_j\bar B_j^*$ has a very simple inverse given by \eqref{eq-UFL-bBbBt}. Thus we can solve the linear system \eqref{eq-M2} efficiently by utilizing the SMW formula described in Section 5.4(a).}

%
%

\begin{center}
	\setlength{\tabcolsep}{1.5pt}
	\tiny 
	\setlength\LTcapwidth{0.9\linewidth}
	\begin{longtable}{| c | c | c | c || c ||  c | c | c | c || c | c | c | c |}
		\hline
		\multicolumn{4}{|c||}{} &\multicolumn{1}{c||}{Original IP} & \multicolumn{4}{c||}{SDP  via DBA-sGS-ADMM} &\multicolumn{4}{c|}{SDP via MOSEK} \\ \hline
		\multicolumn{1}{|c|}{Data} & \multicolumn{1}{c|}{$m_0 \mid m_i$} & \multicolumn{1}{c|}{$n_0 \mid n_i$} & \multicolumn{1}{c||}{$N$}  &
		\multicolumn{1}{c||}{Obj}  & \multicolumn{1}{c|}{Obj} & \multicolumn{1}{c|}{Gap} & \multicolumn{1}{c|}{Iter} &\multicolumn{1}{c||}{Time(s)}  & \multicolumn{1}{c|}{Obj} & \multicolumn{1}{c|}{Gap}  &\multicolumn{1}{c|}{Iter} & \multicolumn{1}{c|}{Time(s)} \\ \hline
		
		\endhead
		
		B1.1  &  51 $\mid$   51 &  51 $\mid$  100 & 100&23468  &23368 &0.43 &1883 &4.08   &23363 &0.45 &  24 &36.33 
		\\[3pt]\hline
		B1.10  &  51 $\mid$   51 &  51 $\mid$  100 & 100&21864  &21193 &3.07 &1911 &4.28   &21187 &3.10 &  24 &36.27 
		\\[3pt]\hline
		C1.1  &  51 $\mid$   51 &  51 $\mid$  100 & 100&16781  &16034 &4.45 &2151 &5.29   &16033 &4.46 &  33 &63.74 
		\\[3pt]\hline
		C1.10  &  51 $\mid$   51 &  51 $\mid$  100 & 100&17994  &17151 &4.68 &1973 &4.40   &17148 &4.70 &  31 &46.04 
		\\[3pt]\hline
		E1.1  &  51 $\mid$   51 &  51 $\mid$  100 & 100&15042  &13876 &7.75 &1969 &3.32   &13871 &7.78 &  28 &44.90 
		\\[3pt]\hline
		E1.10  &  51 $\mid$   51 &  51 $\mid$  100 & 100&14630  &13323 &8.93 &1872 &3.15   &13319 &8.96&  30 &44.99 
		\\[3pt]\hline
		E5.1  &  51 $\mid$   51 &  51 $\mid$  100 & 100&32377  &29290 &9.53 &1603 &2.52   &29279 &9.57 &  25 &37.82
		\\[3pt]\hline
		E5.10  &  51 $\mid$   51 &  51 $\mid$  100 & 100&34086  &30122 &11.63 &1854 &2.84   &30114 &11.65 &  26 &39.18
		\\[3pt]\hline
		E10.1  &  51 $\mid$   51 &  51 $\mid$  100 & 100&46832  &42496 &9.26 &1321 &2.03   &42492 &9.27 &  26 &39.62 
		\\[3pt]\hline
		E10.10  &  51 $\mid$   51 &  51 $\mid$  100 & 100&47449  &43179 &9.00 &1869 &2.96   &43177 &9.00 &  24 &37.82 
		\\[3pt]\hline

		\caption{
			Comparison of computational results for the DNN relaxation of UFL problem in the small scale \texttt{Bilde-Krarup} dataset using DBA-sGS-ADMM and MOSEK. {The results of MOSEK are obtained using \textbf{8 threads}.} Under the column ``Gap'', the relaxation gap is also reported.}
		\label{table:SUFL-SDP}
	\end{longtable}
\end{center}

%
%

\begin{center}
	\setlength{\tabcolsep}{1.5pt}
	\tiny 
	\setlength\LTcapwidth{0.9\linewidth}
	\begin{longtable}{| c | c | c | c || c || c | c | c  |c|}
		\hline
		\multicolumn{4}{|c||}{} &\multicolumn{1}{c||}{Original IP} & \multicolumn{4}{c|}{SDP  via DBA-sGS-ADMM} \\ \hline
		\multicolumn{1}{|c|}{Data} & \multicolumn{1}{c|}{$m_0 \mid m_i$} & \multicolumn{1}{c|}{$n_0 \mid n_i$} & \multicolumn{1}{c||}{$N$}  &
		\multicolumn{1}{c||}{Obj}  & \multicolumn{1}{c|}{Obj} & \multicolumn{1}{c|}{Gap} & \multicolumn{1}{c|}{Iter} &\multicolumn{1}{c|}{Time(s)} 
		\\ \hline
		
		\endhead
		
		ga250a-1  & 251 $\mid$  251 & 251 $\mid$  500 & 250&2.5796e+05 &2.5764e+05 &0.12 &1751 &16.58 
		\\[3pt]\hline
		ga250b-2  & 251 $\mid$  251 & 251 $\mid$  500 & 250&2.7514e+05 &2.7272e+05 &0.88 &1646 &15.61 
		\\[3pt]\hline
		ga250c-3  & 251 $\mid$  251 & 251 $\mid$  500 & 250&3.3366e+05 &3.2288e+05 &3.23 &1504 &13.55  
		\\[3pt]\hline
		ga500a-4  & 501 $\mid$  501 & 501 $\mid$ 1000 & 500&5.1105e+05 &5.1043e+05 &0.12 &2601 &127.43 
		\\[3pt]\hline
		ga500b-5  & 501 $\mid$  501 & 501 $\mid$ 1000 & 500&5.3748e+05 &5.3300e+05 &0.83 &2901 &137.53  
		\\[3pt]\hline
		ga500c-1  & 501 $\mid$  501 & 501 $\mid$ 1000 & 500&6.2136e+05 &6.0286e+05 &2.98 &1796 &83.73 
		\\[3pt]\hline
		ga750a-2  & 751 $\mid$  751 & 751 $\mid$ 1500 & 750&7.6367e+05 &7.6254e+05 &0.15 &3001 &368.26  
		\\[3pt]\hline
		ga750b-3  & 751 $\mid$  751 & 751 $\mid$ 1500 & 750&7.9613e+05 &7.8965e+05 &0.81 &5601 &795.05  
		\\[3pt]\hline
		ga750c-4  & 751 $\mid$  751 & 751 $\mid$ 1500 & 750&9.0004e+05 &8.7556e+05 &2.72 &2197 &246.48 
		\\[3pt]\hline
		gs250a-1  & 251 $\mid$  251 & 251 $\mid$  500 & 250&2.5796e+05 &2.5769e+05 &0.10 &2091 &23.28 
		\\[3pt]\hline
		gs250b-2  & 251 $\mid$  251 & 251 $\mid$  500 & 250&2.7568e+05 &2.7305e+05 &0.95 &1578 &16.92  
		\\[3pt]\hline
		gs250c-3  & 251 $\mid$  251 & 251 $\mid$  500 & 250&3.3300e+05 &3.2199e+05 &3.31 &1707 &17.28 
		\\[3pt]\hline
		gs500a-4  & 501 $\mid$  501 & 501 $\mid$ 1000 & 500&5.1114e+05 &5.1042e+05 &0.14 &3101 &143.43 
		\\[3pt]\hline
		gs500b-5  & 501 $\mid$  501 & 501 $\mid$ 1000 & 500&5.3827e+05 &5.3312e+05 &0.96 &3401 &149.40  
		\\[3pt]\hline
		gs500c-1  & 501 $\mid$  501 & 501 $\mid$ 1000 & 500&6.2004e+05 &6.0198e+05 &2.91 &1871 &81.54 
		\\[3pt]\hline
		gs750a-2  & 751 $\mid$  751 & 751 $\mid$ 1500 & 750&7.6355e+05 &7.6257e+05 &0.13 &3251 &382.80 
		\\[3pt]\hline
		gs750b-3  & 751 $\mid$  751 & 751 $\mid$ 1500 & 750&7.9659e+05 &7.8995e+05 &0.83 &4601 &510.33 
		\\[3pt]\hline
		gs750c-4  & 751 $\mid$  751 & 751 $\mid$ 1500 & 750&9.0134e+05 &8.7541e+05 &2.88 &2608 &292.61 
		\\[3pt]\hline

		\caption{
			Computational results for the DNN relaxation of UFL problem in the large scale \texttt{Koerkel-Ghosh} dataset using DBA-sGS-ADMM. Under the column``Gap'', the relaxation gap is also reported.}
		\label{table:LUFL-SDP}
	\end{longtable}
\end{center}

%

From Table \ref{table:SUFL-SDP}, we could observe that DBA-sGS-ADMM always outperform MOSEK in solving the DNN problems \eqref{eq-UFL-DNN} with dual block angular structure. The runtime is about 10-20 times faster than that  required by MOSEK.

In Table \ref{table:LUFL-SDP}, we report the performance of 
DBA-sGS-ADMM  in solving \eqref{eq-UFL-DNN} 
for the large scale dataset \texttt{Koerkel-Ghosh}.
We can see that our DBA-sGS-ADMM can solve these large scale dual block 
angular DNN problems very efficiently.
In particular, for the last instance {\tt gs750c-4} which has
$564001$ linear constraints, a first stage $751$ by $751$ DNN matrix variables, and 
$750$ second stage $1500$-dimensional nonnegative variables, 
we can solve  it in about 5 minutes. 

We also tried to use MOSEK to solve the instances
in Table \ref{table:LUFL-SDP}, but they often run out of memory on our machine. In addition, due to the inefficiency of MOSEK as observed in Table \ref{table:SUFL-SDP} in solving this type of problems, we do not report the numerical results of MOSEK for these large scale instances.

\subsubsection{DNN relaxations of quadratic UFL problems}

In this subsection, we generate some random quadratic UFL problems following a procedure mentioned in \cite{Gunluk}. 
We compare the performance of our DBA-sGS-ADMM and MOSEK for solving 
the DNN relaxation of the quadratic UFL in \eqref{eq-UFL-DNN} in
Table \ref{table:QUFL-SDP}.
The strategy we adopt to solve the linear system \eqref{eq-M2} is the same as that in Section 7.3.1.

We should take note that MOSEK could not accept the model with both quadratic objective function  and semidefinite constraint. 
Thus we reformulate \eqref{eq-UFL-DNN} as follows.
\blue{First, we have the Cholesky factorization for each ${\rm diag}(Q_j):= H_j^T H_j$, where $H_j={\rm diag}\big(\left(\sqrt{Q_{j1}},...,\sqrt{Q_{jp}}\right)^T\big)$}, then we can reformulate \eqref{eq-UFL-DNN} to be a {convex conic} programming model with rotated quadratic cone constraints:

\begin{eqnarray*}
	\begin{array}{rl}
		\min & \inprod{[0; c]}{[\alp; u]} + \sum_{j=1}^q \inprod{\blue{P_j}}{S_j} + 
		\sum_{j=1}^q t_j 
		\\[8pt]
		\mbox{s.t} &
		\left[\begin{array}{c} 0 \\ u \end{array}\right] 
		+
		\left[\begin{array}{cc}
			e^T &  0^T \\
			-I_p         & -I_p
		\end{array}\right] \left[\begin{array}{c} S_j \\ Z_j \end{array}\right]
		= 
		\left[\begin{array}{c} 1 \\ 0 \end{array}\right], \quad j = 1,\ldots q
		\\[15pt]
		& u - \diag{U} = 0, \;\;  \alp = 1,
		\\[5pt]
		& \mb{U} := \left[\begin{array}{cc} \alp & u^T \\[3pt] u & U \end{array}\right] \in \mathbb{S}^{1+p}_+, \;\;
		\mb{U} \geq 0,\;\; S_j, Z_j \geq 0, \; j=1,\ldots, q,
		\\[15pt]
		&  y_j = H_j S_j, \quad j = 1,\ldots q, 
		\\[5pt]
		&   \beta_j = 1, \quad
		\left[\begin{array}{c}
			\beta_j \\ t_j \\ y_j 
		\end{array}\right] \in \cQ_j:=\{x\in\mathbb{R}^{p+2}: 2x_1 x_2 \ge \sum_{i=3}^{p+2} x_i^2, \; x_1\ge0, \; x_2 \ge 0\}.
	\end{array}
\end{eqnarray*}
We solve the problem \eqref{eq-UFL-DNN} by MOSEK via the above reformulation. 

%
%

\begin{center}
	\setlength{\tabcolsep}{1.5pt}
	\tiny 
	\setlength\LTcapwidth{0.9\linewidth}
	\begin{longtable}{| c | c | c | c ||  c | c | c || c | c | c |}
		\hline
		\multicolumn{4}{|c||}{} &\multicolumn{3}{c||}{SDP  via DBA-sGS-ADMM} &\multicolumn{3}{c|}{SDP  via MOSEK}\\ \hline
		\multicolumn{1}{|c|}{Data} & \multicolumn{1}{c|}{$m_0 \mid m_i$} & \multicolumn{1}{c|}{$n_0 \mid n_i$} & \multicolumn{1}{c||}{$N$}  &
		\multicolumn{1}{c|}{Obj}  & \multicolumn{1}{c|}{Iter} &\multicolumn{1}{c||}{Time(s)} & \multicolumn{1}{c|}{Obj}  & \multicolumn{1}{c|}{Iter} &\multicolumn{1}{c|}{Time(s)}\\ \hline
		
		\endhead
		
		randQUFL-m50-n10  &  11 $\mid$   11 &  11 $\mid$   20 &  50 &83.73  & 101 &0.16    &83.73  &   7 &0.46 
		\\[3pt]\hline
		randQUFL-m50-n20  &  21 $\mid$   21 &  21 $\mid$   40 &  50 &68.11  & 107 &0.10    &68.11  &   10 &0.77 
		\\[3pt]\hline
		randQUFL-m50-n30  &  31 $\mid$   31 &  31 $\mid$   60 &  50 &45.78  & 173 &0.15    &45.78  &  13 &1.83 
		\\[3pt]\hline
		randQUFL-m200-n50  &  51 $\mid$   51 &  51 $\mid$  100 & 200 &73.52  & 402 &0.54    &73.52  &  9 &67.24 
		\\[3pt]\hline
		randQUFL-m200-n100  & 101 $\mid$  101 & 101 $\mid$  200 & 200 &46.97  & 385 &1.01    &46.97  &  16 &1325.55 
		\\[3pt]\hline
		randQUFL-m200-n150  & 151 $\mid$  151 & 151 $\mid$  300 & 200 &42.00  & 469 &1.66   &42.00  &  16 &41258.92 
		\\[3pt]\hline
		randQUFL-m200-n200  & 201 $\mid$  201 & 201 $\mid$  400 & 200 &38.35  & 501 &2.52   & -&- &-
		\\[3pt]\hline
		randQUFL-m1000-n200  & 201 $\mid$  201 & 201 $\mid$  400 &1000 &79.85  & 611 &7.17   &- &- &-
		\\[3pt]\hline
		randQUFL-m1000-n500  & 501 $\mid$  501 & 501 $\mid$ 1000 &1000 &47.94  &1021 &55.13   &- &- &-
		\\[3pt]\hline
		randQUFL-m1000-n1000  &1001 $\mid$ 1001 &1001 $\mid$ 2000 &1000 &35.29  & 791 &159.92   &- &- &-
		\\[3pt]\hline
		randQUFL-m1000-n2000  &2001 $\mid$ 2001 &2001 $\mid$ 4000 &1000 &24.95  & 651 &701.38   &- &- &-
		\\[3pt]\hline

		\caption{
			Comparison of computational results for the DNN relaxation of random large scale QUFL problem using DBA-sGS-ADMM and MOSEK. {The results of MOSEK are obtained using \textbf{8 threads}.} ``-'' means that the solver runs out of time.}
		\label{table:QUFL-SDP}
	\end{longtable}
\end{center}

The comparison of DBA-sGS-ADMM and MOSEK in
Table \ref{table:QUFL-SDP} are similar to that reported in section 7.3.1. One can also observe that DBA-sGS-ADMM is much faster than MOSEK when solving the DNN problems. MOSEK can only solve small or moderate scale problems, and become extremely 
{inefficient} when the scale of problems is large.
{In particular, for the last instance solved by MOSEK, our method is over 24,000 times faster than MOSEK.}

\subsection{Randomly generated problems}
To demonstrate the generality and superiority of our algorithm, in this section, we conduct numerical experiments on some randomly generated dual block angular problems. We have two classes of random problems, namely, the DBA quadratic programming (QP) and the DBA semidefinite programming (SDP).

\subsubsection{QP} 
Here we generate some random quadratic programming problems with dual block angular structure. In particular, we have the following: 
\begin{eqnarray*}
	\theta(x):= \frac{1}{2} \inprod{x}{Qx} \; \text{ and } 
	\btheta_i(\bx_i):= \frac{1}{2} \inprod{\bx_i}{\bQ_i\bx_i} \; \forall i = 1,..., N, \\
	\text{with the bounds } \cK:=\RR^n_+, \; \bcK_i:=\RR^{n_i}_+ \; \forall i = 1,..., N.
\end{eqnarray*}

\noindent 
We randomly generate sparse matrices $A, B$ (by \texttt{MATLAB} command \texttt{sprand}) with density $10/{n_0}$. Similarly, $\bB$ is generated with the same command with density $10/{n_i}$. Finally, we randomly generate sparse symmetric positive definite matrices $Q$ and $\bQ_i$ with density $2/{n_0}$ and $2/{n_i}$ respectively (using \texttt{MATLAB} command \texttt{sprandsym}). {The strategy we adopt to solve the linear system \eqref{eq-M2} is the same as that in Section 7.2.}
\blue{
\begin{center}
	\tiny 
	\setlength{\tabcolsep}{2pt}
	\setlength\LTcapwidth{0.9\linewidth}
	\begin{longtable}{|| c | c | c | c || c | c | c || c | c | c |}
		\hline
		\multicolumn{4}{||c||}{} & \multicolumn{3}{|c||}{DBA-sGS-ADMM} & \multicolumn{3}{|c|}{Gurobi} \\ \hline
		\multicolumn{1}{||c|}{Data} & \multicolumn{1}{|c|}{$m_0 \mid m_i$} & \multicolumn{1}{|c|}{$n_0 \mid n_i$}
		& \multicolumn{1}{|c||}{$N$}  & \multicolumn{1}{|c|}{Iter} &\multicolumn{1}{|c|}{Time(s)} & \multicolumn{1}{|c||}{\objP} 	&\multicolumn{1}{|c|}{Iter} &\multicolumn{1}{|c|}{Time(s)} & \multicolumn{1}{|c|}{\objP}  \\ \hline
		
		\endhead
		
		randQP-m10-n20-N10  &  10 $\mid$   10 &  20 $\mid$   20 &  10 & 104   &0.11 &3.523e+02 &  11 &0.28 &3.523e+02
		\\[3pt]\hline
		randQP-m50-n80-N10  &  50 $\mid$   50 &  80 $\mid$   80 &  10 & 151   &0.15 &1.806e+03 &  10 &0.30 &1.806e+03
		\\[3pt]\hline
		randQP-m100-n200-N10  & 100 $\mid$  100 & 200 $\mid$  200 &  10 & 283   &0.50 &3.707e+03 &  11 &0.51 &3.707e+03
		\\[3pt]\hline
		randQP-m200-n300-N10  & 200 $\mid$  200 & 300 $\mid$  300 &  10 & 374   &2.77 &7.048e+03 &  11 &1.66 &7.048e+03
		\\[3pt]\hline
		randQP-m500-n800-N10  & 500 $\mid$  500 & 800 $\mid$  800 &  10 & 335   &3.77 &1.770e+04 &  10 &17.34 &1.770e+04
		\\[3pt]\hline
		randQP-m100-n200-N50  & 100 $\mid$  100 & 200 $\mid$  200 &  50 & 810   &2.50 &1.713e+04 &  13 &0.99 &1.713e+04
		\\[3pt]\hline
		randQP-m200-n300-N50  & 200 $\mid$  200 & 300 $\mid$  300 &  50 & 957   &7.02 &3.362e+04 &  10 &3.66 &3.362e+04
		\\[3pt]\hline
		randQP-m500-n800-N50  & 500 $\mid$  500 & 800 $\mid$  800 &  50 & 791   &32.86 &8.494e+04 &  11 &20.23 &8.494e+04
		\\[3pt]\hline
		randQP-m100-n200-N100  & 100 $\mid$  100 & 200 $\mid$  200 & 100 & 973   &5.65 &3.460e+04 &  13 &7.45 &3.460e+04
		\\[3pt]\hline
		randQP-m200-n300-N100  & 200 $\mid$  200 & 300 $\mid$  300 & 100 &1023   &14.57 &6.687e+04 &  12 &11.53 &6.687e+04
		\\[3pt]\hline
		randQP-m500-n800-N100  & 500 $\mid$  500 & 800 $\mid$  800 & 100 & 601   &7.30 &1.688e+05 &  13 &70.68 &1.688e+05
		\\[3pt]\hline
		randQP-m100-n200-N200  & 100 $\mid$  100 & 200 $\mid$  200 & 200 &1169   &13.35 &6.897e+04 &  15 &3.41 &6.897e+04
		\\[3pt]\hline
		randQP-m200-n300-N200  & 200 $\mid$  200 & 300 $\mid$  300 & 200 &1187   &33.57 &1.324e+05 &  13 &19.67 &1.324e+05
		\\[3pt]\hline
		randQP-m500-n800-N200  & 500 $\mid$  500 & 800 $\mid$  800 & 200 & 859   &18.58 &3.369e+05 &  13 &118.72 &3.369e+05
		\\[3pt]\hline
		randQP-m1000-n2000-N200  &1000 $\mid$ 1000 &2000 $\mid$ 2000 & 200 &1326   &68.79 &7.003e+05 &  22 &869.92 &7.003e+05
		\\[3pt]\hline
		randQP-m2000-n3000-N200  &2000 $\mid$ 2000 &3000 $\mid$ 3000 & 200 & 916  &114.02 &1.351e+06 &  16 &4237.40 &1.351e+06
		\\[3pt]\hline
		randQP-m5000-n8000-N200  &5000 $\mid$ 5000 &8000 $\mid$ 8000 & 200 & 901   &466.62 &3.403e+06&- & -&-
		\\[3pt]\hline

		\caption{
			Comparison of computational results between DBA-sGS-ADMM and Gurobi for randomly generated \textbf{QP} problem. All the results are obtained using a {\textbf{8 thread}}. The objective function value is shown to help us verifying the correctness easily. ``-'' means that the solver runs out of time.}
		\label{table:randQP}
	\end{longtable}
\end{center}
}

From Table \ref{table:randQP}, we can observe that DBA-sGS-ADMM is much faster than Gurobi when $m_0,m_i,n_0,n_i$ are large. We note that $n_0,n_i$ represent the dimensions of first stage and second stage decision variables respectively. Meanwhile, we generate denser $Q$ in the quadratic objective than that in section 7.2.2 ($Q=0.1I$), making the computational cost of Gurobi per iteration become more expensive. Therefore, we could expect that our DBA-sGS-ADMM 
to be more efficient in solving two stage stochastic quadratic programming problems when 
the matrix in quadratic objective is dense and the dimensions of the decision variables are large.

\subsubsection{SDP}
We also generate a random data set of semidefinite programming problems with dual block angular structure. In particular, we have
\begin{eqnarray*}
	&&\theta(x):= 0, \;
	\btheta_i(\bx_i):=0 \; \forall i = 1,..., N, \\
	&&\text{with the bounds } \cK:=S^n_+, \; \bcK_i:=S^{n_i}_+ \; \forall i = 1,..., N, \\
	&&\text{and the mappings }
	A(x)
	:=\begin{bmatrix}
		\inprod{A_1}{x} \\ \inprod{A_2}{x} \\ \vdots \\ \inprod{A_{m}}{x}
	\end{bmatrix}, \;
	B_i(x)
	:=\begin{bmatrix}
		\inprod{B_{i,1}}{x} \\ \inprod{B_{i,2}}{x} \\ \vdots \\ \inprod{B_{i,m_i}}{x}
	\end{bmatrix},  \\
	&&\bar{B}_i(\bx_i)
	:=\begin{bmatrix}
		\inprod{\bar{B}_{i,1}}{\bx_i} \\
		\inprod{\bar{B}_{i,2}}{\bx_i} \\
		\vdots \\
		\inprod{\bar{B}_{i,m_i}}{\bx_i} 
	\end{bmatrix}.
\end{eqnarray*}	

\noindent 
We generate the matrix representation of the linear mapping $A_i$ for $i=1,\ldots,m$ (using \texttt{MATLAB} routine 
{\texttt{Ai = sprand(n0,n0,0.2); Ai = Ai*Ai';}}).
Similarly, the matrix representation of $B_{i,k}$ and $\bar{B}_{i,k}$ for  $k=1,\ldots,{m}_i$,
$i=1,\ldots,N$, is generated using the same routine except with density $5/{n_i}$. {The strategy we adopt to solve the linear system \eqref{eq-M2} is the same as that in Section 7.2.}

\begin{center}

	\tiny 
	\setlength{\tabcolsep}{2pt}
	\setlength\LTcapwidth{0.9\linewidth}
	\begin{longtable}{|| c | c | c | c ||  c | c | c || c | c | c |}
		\hline
		\multicolumn{4}{||c||}{} & \multicolumn{3}{|c||}{DBA-sGS-ADMM} & \multicolumn{3}{|c|}{MOSEK} \\ \hline
		\multicolumn{1}{||c|}{Data} & \multicolumn{1}{|c|}{$m_0 \mid m_i$} & \multicolumn{1}{|c|}{$n_0 \mid n_i$}
		& \multicolumn{1}{|c||}{$N$}  & \multicolumn{1}{|c|}{Iter} &\multicolumn{1}{|c|}{Time(s)} & \multicolumn{1}{|c||}{\objP} 	&\multicolumn{1}{|c|}{Iter} &\multicolumn{1}{|c|}{Time(s)} & \multicolumn{1}{|c|}{\objP}  \\ \hline
		
		\endhead
		
		randSDP-m50-n50-N5  &  50 $\mid$   50 &  100 $\mid$   100 &  5 &544   &8.79 &2.529e+05 &   10 &2.72 &2.529e+05
		\\[3pt]\hline
		randSDP-m100-n100-N5  &  100 $\mid$   100 &  100 $\mid$   100 &  5 &202   &5.56 &4.682e+05 &   9 &6.89 &4.682e+05
		\\[3pt]\hline
		randSDP-m100-n200-N5  &  100 $\mid$   100 &  200 $\mid$   200 &  5 &501   &36.22 &2.769e+06 &   13 &32.77 &4.778e+06
		\\[3pt]\hline
		randSDP-m200-n200-N5  &  200 $\mid$   200 &  200 $\mid$   200 &  5 &201   &30.35 &4.778e+06 &   10 &57.26 &4.778e+06
		\\[3pt]\hline
		randSDP-m200-n300-N5  &  200 $\mid$   200 &  300 $\mid$   300 &  5 &256   &70.95 &1.521e+07 &   11 &159.41 &1.521e+07
		\\[3pt]\hline
		randSDP-m200-n400-N5  &  200 $\mid$   200 &  400 $\mid$   400 &  5 &277   &125.52 &3.710e+07 &   13 &389.44 &3.710e+07
		\\[3pt]\hline
		randSDP-m400-n400-N5  &  400 $\mid$   400 &  400 $\mid$   400 &  5 &265   &248.69 &6.399e+07 &   10 &838.80 &6.398e+07
		\\[3pt]\hline
		randSDP-m10-n20-N10  &  10 $\mid$   10 &  20 $\mid$   20 &  10 &2375   &3.48 &5.222e+03 &   9 &0.49 &5.222e+03
		\\[3pt]\hline
		randSDP-m50-n80-N10  &  50 $\mid$   50 &  80 $\mid$   80 &  10 & 236   &4.88 &2.327e+05 &  10 &3.37 &2.327e+05
		\\[3pt]\hline
		randSDP-m50-n100-N10  &  50 $\mid$   50 &  100 $\mid$   100 &  10 & 162   &4.80 &3.631e+05 &  11 &5.58 &3.631e+05
		\\[3pt]\hline
		randSDP-m80-n100-N10  &  80 $\mid$   80 &  100 $\mid$   100 &  10 & 240   &9.67 &5.406e+05 &  10 &10.69 &5.406e+05
		\\[3pt]\hline
		randSDP-m100-n200-N10  & 100 $\mid$  100 & 200 $\mid$  200 &  10 & 212   &28.73 &3.286e+06 &  11 &48.13 &3.286e+06
		\\[3pt]\hline
		randSDP-m200-n300-N10  & 200 $\mid$  200 & 300 $\mid$  300 &  10 & 309   &128.05 &1.646e+07 &  13 &379.36 &1.646e+07
		\\[3pt]\hline
		randSDP-m200-n400-N10  & 200 $\mid$  200 & 400 $\mid$  400 &  10 & 280   &190.98 &3.907e+07 &  12 &773.04 &3.907e+07
		\\[3pt]\hline
		randSDP-m200-n400-N20  & 200 $\mid$  200 & 400 $\mid$  400 &  20 & 633   &680.14 &3.906e+07 &  12 &2042.76 &3.906e+07
		\\[3pt]\hline 

		\caption{
			Comparison of computational results between DBA-sGS-ADMM and MOSEK for randomly generated \textbf{SDP} problem. The results of MOSEK are obtained using \textbf{8 threads}. }
		\label{table:randSDP}
	\end{longtable}

\end{center}

In Table \ref{table:randSDP}, the number of scenarios $N$ is small. When the number of scenarios $N$ is fixed and the dimensions of the decision variables increase gradually, we can observe that the results are consistent with the SDP relaxation problems in section 7.3. 
When $m_0,m_i,n_0,n_i$ increase further, our DBA-sGS-ADMM 
is expected to outperform MOSEK by an even larger margin.

\section{Conclusion}

{
	We have proposed efficient augmented Lagrangian based decomposition algorithms for solving 
	convex composite programming problems with dual block-angular structures. 
	We show that various application problems such as 2-stage stochastic convex programming problems, and DNN relaxations of uncapacitated facility location problems can be recast into our framework.  Numerical experiments show that our method is especially competitive against MOSEK or Gurobi
	for solving SDP problems and
	convex QP problems with dual block-angular structures. Our algorithms are also 
	demonstrated to be much more competitive than the well-known PHA in solving 
	2-stage stochastic linear programming problems.
	As mentioned in the paper, some of the subproblems can be solved in a parallel manner. So as a natural future work, we may implement our algorithm on
	a conducive parallel computing environment to further improve the efficiency of our algorithms.
}

\section*{Acknowledgements}
We would like to thank the two anonymous referees for their careful reading of this
paper, and their insightful comments and suggestions which have helped to improve the quality of this paper.

\begin{appendices}
	
	\section*{Appendix 1: two-stage stochastic programming datasets }
	\begin{description}
		\item[assets.] This network model represents the management of assets. Its nodes are asset categories and its arcs are transactions. The problem is to maximize the return of an investment from every stage with the balance of material at each node.
		
		\item[env.] This model assists the Canton of Geneva in planning its energy	supply infrastructure and policies. The main objective is to minimize the installation cost of various types of energy, while meeting all the supply-demand at each node and satisfying several realistic constraints such as the environmental constraints. 
		
		\item[phone.] This is a problem which models the service of providing private lines to telecommunication customers, often used by large corporations between business locations for high speed, private data transmission. The goal is to minimize the unserved requests,	while staying within budget.
		
		\item[AIRL.] This dataset is used to schedule monthly airlift operations. The airlift operation is scheduled so that the restriction on each aircraft such as the number of flight hours available during the month can be satisfied. The cost such as the available flight time to go unused, switching aircraft from one route to another and buying commercial flights is the main objective to be minimized.
		
		\item[4node.] This is a two stage network problem for scheduling cargo transportation. While the flight schedule is completely determined in stage one, the amounts of cargo to be shipped are uncertain and shall be determined in stage two. 
		
		\item[pltexp.] This is a stochastic capacity expansion model that tries to allocate new production capacity across a set of plants so as to maximize profit subject to uncertain demand. 
		
		\item[storm.] This is a two period freight scheduling problem described in Mulvey and Ruszczynski. In this model, routes are scheduled to satisfy a set of demands at stage 1, demands occur, and unmet demands are delivered at higher costs in stage 2 to account for shortcomings.
		
		\item[gbd.] This is the aircraft allocation problem where aircraft of different types are to be allocated to routes in a way that maximizes profit under uncertain demand, and minimizes the cost of operating the aircraft as well as costs associated with bumping passengers when the demand for seats outstrips the capacity. 
	\end{description}

	\section*{Appendix 2: Derivation of solutions of \eqref{ALMSNN-subp2-zy} and \eqref{ALMSNN-subp2-yv}}
\blue{
We can rewrite the subproblem in \eqref{ALMSNN-subp2-zy} as follows:
	\begin{eqnarray*}
	\begin{aligned} 
	&\min\limits_{z,y}  \Big\{ -\inprod{b}{y} + \delta_\cK^*(-z) + 
	\frac{\sigma}{2}\norm{z+A^*y - \widehat{c}^k}^2 -\inprod{\delta^{k}}{y}
	\Big\}\\
	=&\min\limits_{y}\left\{-\inprod{b}{y}-\inprod{\delta^k}{y}+\min\limits_{z}\left\{\delta^*_{\cK}(-z)+\frac{\sigma}{2}\|z+A^*y-\hat c^k\|^2\right\}\right\}\\
	=&\min\limits_{y}\left\{-\inprod{b}{y}-\inprod{\delta^k}{y}+\sigma M_{\delta^*_{\cK}/\sigma}(A^*y-\hat c^k)\right\},
	\end{aligned}
	\quad
	\end{eqnarray*}
	where the last equality comes from the definition of Moreau-Yosida envelope, and the equality holds for $z=-\mathop{\rm{Prox}}_{\delta^*_{\cK}/\sigma}(A^*y-\hat c^k)$. By the property of Moreau-Yosida proximal mapping, $\mathop{\rm{Prox}}_{\delta^*_{\cK}/\sigma}(A^*y-\hat c^k)=(A^*y-\hat c^k)-\mathop{\rm{Prox}}_{\sigma\delta_{\cK}}(\sigma(A^*y-\hat c^k))$, therefore, $z=\Pi_{\cK}(\sigma(A^*y-\hat c^k))-(A^*y-\hat c^k)$, which gives \eqref{ALMSNN-subp2-zy}. The derivation of \eqref{ALMSNN-subp2-yv} is similar, we omit the details.
	
}

\end{appendices}




\subsection*{Statements and Declarations}
The research of Kim-Chuan Toh is supported
by the Ministry of Education, Singapore, under its Academic Research Fund Tier 3 grant call (MOE-2019-T3-1-010). All authors certify that they have no affiliations with or involvement in any organization or entity with any financial interest or non-financial interest in the subject matter or materials discussed in this manuscript.

\subsection*{Data availability statement}
The links of datasets from the literature that are used in this paper have been included where they are mentioned. For the other simulated data, they are available from the corresponding author on reasonable request.

\end{document}